\documentclass{article}
\usepackage[utf8]{inputenc}
\usepackage[left=3cm, right=3cm,top=3cm,bottom=3cm]{geometry}
\usepackage{amsmath}
\usepackage{amssymb}
\usepackage{amsfonts}
\usepackage{amsthm}
\usepackage{tikz-cd}
\usepackage{calligra,mathrsfs}
\usepackage{tikz-cd}
\usepackage{stmaryrd}
\usepackage{cleveref,mathtools}

\newcommand{\comment}[1]{}

\theoremstyle{plain} 
\newtheorem{thm}{Theorem}[section]
\newtheorem{prop}[thm]{Proposition}
\newtheorem{lem}[thm]{Lemma}
\newtheorem{cor}[thm]{Corollary}
\newtheorem*{introcor}{Corollary}

\newtheorem{hypo}[thm]{Assumption}
\newtheorem*{thm*}{Theorem}
\newtheorem{introthm}{Theorem}[section]

\theoremstyle{definition}  

\newtheorem{defn}[thm]{Definition}
\newtheorem{notation}[thm]{Notation}

\newtheorem*{rem*}{Remark}

\theoremstyle{remark}
\newtheorem{rem}[thm]{Remark}

\newcommand{\N}{\mathbb{N}}
\newcommand{\Z}{\mathbb{Z}}
\newcommand{\Q}{\mathbb{Q}}

\newcommand{\C}{\mathbb{C}}
\newcommand{\A}{\mathbb{A}}
\newcommand{\G}{\mathbb{G}}

\newcommand{\td}{\textrm{td}}

\newcommand{\GW}{\mathrm{GW}}

\newcommand{\spine}{\mathrm{spine}}
\newcommand{\tail}{\mathrm{tail}}
\newcommand{\un}{\mathrm{un}}
\newcommand{\pr}{\mathrm{pr}}

\newcommand{\FJRW}{\mathrm{FJRW}}
\newcommand{\pt}{\mathrm{pt}} 

\DeclareMathOperator{\Res}{Res}

\allowdisplaybreaks[3]

\newcommand{\arm}{\mathrm{arm}}
\newcommand{\leg}{\mathrm{leg}}
\newcommand{\head}{\mathrm{head}}

\newcommand{\lcm}{\mathrm{lcm}}

\newcommand{\ua}{\underline{a}}
\newcommand{\ub}{\underline{b}}

\newcommand{\ug}{\underline{\gamma}}

\newcommand{\fa}{\mathfrak{a}}

\newcommand{\fe}{\mathfrak{e}}

\newcommand{\PP}{\mathbb{P}}

\newcommand{\cO}{\mathcal{O}}
\newcommand{\ev}{\mathrm{ev}}

\newcommand{\cA}{\mathcal{A}}
\newcommand{\cB}{\mathcal{B}}
\newcommand{\cC}{\mathcal{C}}
\newcommand{\cD}{\mathcal{D}}

\newcommand{\cF}{\mathcal{F}}
\newcommand{\cH}{\mathcal{H}}
\newcommand{\cI}{\mathcal{I}}
\newcommand{\scrI}{\mathscr{I}}
\newcommand{\cJ}{\mathcal{J}}

\newcommand{\cK}{\mathcal{K}}

\newcommand{\cL}{\mathcal{L}}
\newcommand{\cM}{\mathcal{M}}
\newcommand{\cN}{\mathcal{N}}

\newcommand{\cT}{\mathcal{T}}

\newcommand{\cX}{\mathcal{X}}
\newcommand{\cY}{\mathcal{Y}}

\newcommand{\bmu}{\boldsymbol{\mu}}

\newcommand{\Td}{\mathrm{Td}}

\newcommand{\Spec}{\mathrm{Spec}}

\DeclareMathOperator{\ch}{ch}

\newcommand{\ff}{\boldsymbol{f}}

\DeclareMathOperator{\Aut}{Aut}

\DeclareMathOperator{\cHom}{\mathscr{H}\text{\kern -3pt {\calligra\large om}}\,}

\newcommand{\tr}{\operatorname{tr}}
\newcommand{\Tr}{\operatorname{Tr}}
\newcommand{\mult}{\operatorname{mult}}

\def\corr#1{\left\langle #1 \right\rangle}
\newcommand{\scal}[1]{\left\langle #1 \right\rangle}

\def\<{\left\langle}
\def\>{\right\rangle}

\newcommand{\fake}{\mathrm{fake}}

\makeindex
\title{Permutation-equivariant quantum K-theory of Fermat singularities}
\author{Maxime Cazaux}

\begin{document}
\maketitle
\begin{abstract}
We compute the genus-0 permutation-equivariant quantum K-theory of Fermat singularities, in parallel with the Givental--Lee theory for projective varieties.
We extend Givental--Tonita's formalism of adelic Lagrangian cones to the singularity theory, and we obtain explicit $I$-functions for the invariants, which satisfy the same $q$-difference equation as Givental's $I$-function of the associated hypersurface.
This can be regarded as an extension of the Landau--Ginzburg/Calabi--Yau correspondence, although a discrepancy between the two sides sides emerges in K-theory.
In the case of the quintic threefold, both generating functions satisfy a $q$-difference equation of degree $25$; the hypersurface $I$-function only spans a $5$-dimensional subspace of solutions, while the singularity $I$-function spans the full space of solutions. 
\end{abstract}


\section*{Introduction}
The quintic threefold occupies a central place in Gromov--Witten (GW) theory since its very definition \cite{candelasPairCalabiYauManifolds1991}.
Still, even if major mirror symmetry conjectures about it have been proven \cite{giventalEquivariantGromovWittenInvariants1996,lianMirrorPrinciple1999}, and the genus $1$ and $2$ have been successfully treated \cite{zingerReducedGenus$1$2008,guoMirrorTheoremGenus2017}, the full Gromov--Witten potential remains unknown. 
Several methods to compute GW invariants have been found \cite{maulikTopologicalViewGromov2006,changEffectiveTheoryGW2022,arguzGromovWittenTheory2023}, but it is still impossible to reach the conjectural formulae by Huang--Klemm--Quackenbush for $g\leq 51$ \cite{huangTopologicalStringTheory2009}.
Remarkably, we even lack  conjectures  beyond $g=51$.

In order to overcome this difficulty, following ideas of Witten \cite{wittenPhases$N2$Theories1993}, an alternative approach has been extensively developed in the last decades. 
In simple terms, instead of focusing on the hypersurface $X$ within the complex projective space $\PP^4$, we regard the associated affine cone yielding a singularity at the origin of $\A^5$. 
The quantum theory of this singularity, referred to in physics literature as the Landau--Ginzburg model, is the so-called FJRW cohomological field theory (CohFT), constructed by Fan, Jarvis and Ruan in the analytic category \cite{fanWittenEquationMirror2013}, and by Polischchuk and Vaintrob in algebro-geometric terms \cite{polishchukMatrixFactorizationsCohomological2016}.
As Witten explains, this new point of view arises from a change of stability condition in geometric invariant theory for the action of a reductive group $G$ on a vector space enriched with a $G$-equivariant complex-valued function $W$.
This is a general setup referred to as the gauged linear sigma model (GLSM), which in principle recovers the geometry of the quintic hypersurface in $\PP^4$ and that of the singularity in $\A^5$ via a change of stability condition.
Therefore, it is natural to expect the so-called Landau--Ginzburg/Calabi--Yau (LG/CY) correspondence between the FJRW invariants of the singularity (LG side) and the GW invariants of the quintic (CY side). 
Interestingly, the two theories are radically different in their moduli spaces and in the cohomology classes involved; 
therefore, they are expected to shed new light to each other.
This idea found confirmation in a series of results, among which we can 
mention a few that will play a role in the development illustrated here.
Chiodo, Ruan and Iritani \cite{chiodoLandauGinzburgCalabiYauCorrespondence2010,chiodoLandauGinzburgCalabiYauCorrespondence2014} cast it within the framework of the LG/CY 
correspondence, compatibly with Orlov equivalence.
Fan--Jarvis--Ruan \cite{fanMathematicalTheoryGauged2018} constructed a mathematical theory of the GLSM, and Chan--Li--Li--Liu \cite{changEffectiveTheoryGW2022} provided an algorithm computing both GW and FJRW invariants.

The goal of this paper is to extend these methods to quantum K-theory, an analogue of GW theory introduced by A.~Givental and Y.-P.~Lee \cite{leeQuantumKtheoryFoundations2004}.
It turns out that on the singularity side the relevant 
K-theoretic invariant are already defined: the definition of 
the virtual class 
in FJRW theory factors through K-theory 
since its early definition. 
The first example was the 
theory of the $A_{r-1}$ singularity which was 
identified to  the theory of $r$th roots $\cL$ of the (log) canonical bundle;
there, 
Polishchuk--Vaintrob \cite{polishchukAlgebraicConstructionWittens2001} and Chiodo \cite{chiodoWittenTopChern2006a}
provided a definition of the relevant intersection numbers directly in $K$-theory. 
A later construction even produces objects in the derived category of the relevant moduli spaces \cite{polishchukMatrixFactorizationsCohomological2016}. 

We consider the case of the so-called Fermat polynomial $W=\sum_{i=0}^{N}X_{i}^{r}$.
Under the concavity assumption, the K-class for the FJRW invariants of $W$ boils down to the K-theoretic Euler class $\lambda_{-1}$ of the vector bundle 
$R^1\pi_*\cL^{\oplus N+1}$,
where $\pi\colon \cC\to \cM$ is the universal (twisted) curve on the 
moduli space $\cM$ of $r$th roots $\cL$ of the (log) canonical bundle \cite{polishchukMatrixFactorizationsCohomological2016,tsengWallCrossingGenusZero2016,guereCongruencesKtheoreticGromov2023}.  
In this paper, we deal with a refinement of quantum K-theory 
introduced by Givental under the name of  
\emph{permutation-equivariant quantum K-theory} \cite{giventalPermutationEquivariantQuantumKtheory}: 
its invariants encode not only the Euler characteristics, but rather the full $S_n$-module structure of the virtual fundamental sheaf (see \Cref{defn:invariants}). 

For the permutation-equivariant K-theoretic invariants we provide a full computation
encoded in a generating function $I^K_{\FJRW}$ (see \Cref{thm:fonctionI}).
It turns out that, up to a prefactor and a change of variable, $I^{K}_{\FJRW}$ is a solution to the $q$-difference equation satisfied by the generating function $I^K_{\GW}$ of the quantum K-theoretic invariants of the hypersurface $X=\{W=0\}\subset \PP^{N}$ (see \cite{giventalPermutationEquivariantQuantumKtheoryb})
\begin{equation}
\label{intro:qdiffequ0}
\left[(1-q^{Q\partial_Q})^{N+1} -Q\prod_{k=1}^r (1-q^{rQ\partial_Q +k})\right]P^{l_q(Q)}I^K_{\GW}(q,Q) =0.
\end{equation}

In the case of the quintic polynomial (i.e. $r=N+1=5$), $I^{K}_{\FJRW}$ recovers the full space of solutions to \eqref{intro:qdiffequ0}.
This is  an interesting improvement over the Gromov--Witten theory, where the generating function $I^K_{\GW}$ only provides a 5-dimensional subspace of solutions of the degree-25 equation \eqref{intro:qdiffequ0}. 
Here, $25$-dimensions are spanned, and 
this is made possible by the fact that the FJRW invariants are naturally defined over the $r^{2}$-dimensional state space $\C[\widehat{\bmu}_r]^r$; \emph{i.e.}, the K-theory of the inertia stack of the target $B\bmu_r$.
Furthermore, these functions match the basis of solutions already 
identified by Yaoxiong Wen via an analytic continuation of $I_{\GW}^K$ \cite{wenDifferenceEquationQuintic2022}.

A large scale picture beyond the case of the quintic incorporating permutation-equivariant quantum $K$-theory, 
mirror symmetry and the
LG/CY correspondence 
 is still lacking, but it is also 
getting more and more precise thanks to a number of recent papers
over the past ten years. 
We refer to the work of Konstantin Aleshkin and Melissa Liu \cite{aleshkinWallcrossingKtheoreticQuasimap2022,aleshkinHiggsCoulombCorrespondenceWallCrossing2023} 
and references therein. 
There, the authors 
used the above mentioned framework of GLSM 
in order to deduce  
$I$-functions satisfying  $q$-difference equations.
Under some Calabi--Yau conditions they obtain 
wall crossing statements \cite[Thm.~4.3, Prop.~4.6]{aleshkinWallcrossingKtheoreticQuasimap2022} similar to the LG/CY correspondence.   
One can hope that this approach and the present 
paper can contribute to recast the permutation-equivariant K-theory in a 
global mirror symmetry framework as it happens for Gromov--Witten theory in 
\cite{chiodoLandauGinzburgCalabiYauCorrespondence2014}.

\paragraph{Overview of the main results.}
We set up the theory for a Fermat hypersurface $X$ in $\PP^{N}$.
The K-theoretic FJRW invariants of the polynomial $W=X_0^r+\cdots + X_N^r$ are defined via the moduli space $\widetilde{\cM}^r_{0,n}$ of $r$-spin curves with trivialised marking, which parametrizes twisted curves together with an $r$th root $\cL$ of the log-canonical bundle $\omega_{\log}$, and $n$ sections $\sigma_{i}$.
The universal curve $\pi : \widetilde{\cC}^r_{0,n} \to \widetilde{\cM}^r_{0,n}$ carries a universal $r$th root $\cL$, which is used to define evaluation maps $\ev_{i}: \widetilde{\cM}^{r}_{0,n} \to \cI B\bmu_{r}$.
We define the virtual class as the K-theoretic Euler class of the higher direct image of $\cL$, twisted by a certain divisor $E$
\begin{equation}
\label{intro:equ:lambda}
\Lambda_{n} = \lambda_{-1}\left( R^1\pi_*\cL(-E)^{\oplus N+1}\right) \in K^0\left(\widetilde{\cM}^r_{0,n}\right).
\end{equation} 
For K-theoretic classes $F,F_{1},\hdots, F_{n} \in K^{0}\left(\cI B\bmu_{r}\right)$, the K-theoretic FJRW invariants are defined as the Euler characteristics
\begin{equation*}
	\chi \left(\overline{\cM}^r_{0,n}; \Lambda_{n} \bigotimes_{i=1}^{n}\ev_{i}^{*}F_{i}\right),
\end{equation*} 
or, in their permutation-equivariant version, as the $S_n$-module $$H^* \left(\overline{\cM}^r_{0,n}; \Lambda_{n} \bigotimes \ev_{i}^{*}F\right).$$

The computation of these invariants is parallel to Givental--Tonita's computation of the genus-0 quantum K-theory \cite{giventalHirzebruchRiemannRoch2011,giventalPermutationEquivariantQuantumKtheorya} via Lefschetz trace formula, which is an instance of the Grothendieck--Riemann--Roch theorem for stacks \cite{kawasakiRiemannRochTheoremComplex1979,toenTheoremesRiemannRoch1999}.
If $g$ is a finite-order automorphism of a smooth proper stack $X$, and $F$ is an equivariant coherent sheaf, then the trace of $g$ on the cohomology groups of $F$ is given by
$$
\label{intro:equ:lefschetz}
\tr_g(H^*(X, F))= \int_{X^g} \ch\left( \frac{\Tr(F)}{\Tr(\lambda_{-1}\cN^\vee)} \right)\td(T),
$$
where $X^g$ is the fixed-point stack, $\cN$ is the normal bundle to the map $X^g\to X$, $\Tr(F)$ is the trace bundle, and $\td(T)$ is the Todd class of tangent sheaf.
This formula is used to compute all the K-FJRW invariants by recursion on the number of markings.
Indeed, the integral above takes place on the fixed-point stack of $\overline{\cM}^r_{0,n}$, which is the disjoint union of $r$ components of dimension $\dim \overline{\cM}_{0,n}=n-3$, and other lower-dimensional boundary strata with the usual Deligne--Mumford recursive structure.
Such a recursive structure is the core of Givental's formalism where the invariants take the form of a generating function in a polarized symplectic space.
Thus, the full computation derives from inserting the top-dimensional classes $\ch(\lambda_{-\xi}R^1\pi_*(\cL^{\oplus N+1}))$ in Givental's formalism.
These classes are usually referred to as the fake quantum K-theory, and come here in $r$ variants indexed by $\xi \in \mu_r$.
Givental's symplectic space is $\cK = K^0(\cI B\bmu_r)(q^{1/r})$, and the generating function $J:\cK_+ \to \cK$ is defined by
$$J(t)= 1-q + t + \sum_{a,\xi}\sum_{n\geq 2} \phi^{a,\xi} \scal{\frac{\phi_{a,\xi }}{1-q^{1/\fe(a)}\cL_0}, t(\cL_1),\hdots, t(\cL_n)}^{S_n}_{0,n+1}.$$
The state space also carries a natural $\bmu_r$-action, and we obtain an analogue of Givental--Tonita's adelic characterization theorem \cite{giventalHirzebruchRiemannRoch2011}.
\begin{introthm}
Let $f$ be a $\bmu_r$-invariant element of $\cK$.
Then $f$ lies in the image of the $J$-function if and only if
\begin{itemize}
\item $f$ has poles only at $q=0,\infty$, and at the roots of unity;
\item the expansion of $f$ at $q=1$ is a value of the fake $J$-function;
\item for all $\xi_0\in \bmu_\infty$ such that $\xi_0^r$ has order $m$,
	we have $$\Phi_{\xi_0} \left(f(q^{\frac{1}{rm}})\right) \in \Box_{\xi_0}\Delta^{-1}\cT L^\fake $$ where $\cT L^\fake$ is the tangent space at the point $\Phi_0f \in L^\fake$, $\Delta$ is the operator of the fake theory (\Cref{prop:Delta}), $\Box_{\xi_0}$ is the operator of the spine theory (see \Cref{section:spinecohft}), and the maps $\Phi_{\zeta}$ are defined in \Cref{defn:morphismePhi}.

\end{itemize}
\end{introthm}
We deduce a simpler characterization of the $J$-function in terms of the so-called  untwisted invariants.
The untwisted invariants encode the $S_{n}$ modules obtained as cohomology groups of tautological line bundles on the moduli space of $r$-pin curves (see \Cref{subsection:untwisted}).
They are easily determined by the quantum K-theory of the point \cite{giventalPermutationEquivariantQuantumKtheory}.
\begin{introcor}
Let $L^K_{\FJRW}$ be the image of the $J$-function, and let $L^{K}_{\un}$ be the image of the untwisted $J$-function.
Then we have
$$(L)^{\bmu_r}=\Delta (L^{K}_{\un})^{\bmu_r}.$$
\end{introcor}

We use the previous result to find a specific point of $\cL^{K}_{\FJRW}$.
\begin{introthm}
The following function lies on  $L_{\FJRW}^K$ 
$$ I_{\FJRW}^K(x,q)=(1-q)\sum_{\xi\in \bmu_r}  \sum_{n \geq 0 }  \frac{\prod_{0\leq k < n} \left(1-s\xi q^{\{\frac{a+1}{r}\}+k}\right)^{N+1}}{\prod_{k=1}^{n}(1-q^k)} x^{n} \phi_{n+1,\xi}. $$ 
\end{introthm}
We decompose the $I$-function into $r^{2}$ functions
$$I_{\FJRW}^K(x,q) = \sum_{\substack{a=0,\hdots, r-1\\ \xi\in \bmu_{r}}}  x^a I_{a,\xi}(x,q) \phi_{a+1,\xi}.$$
We introduce the modification $$\tilde{I}_{a,\xi}(x,q^{-1})=e_{q,q^{\frac{a+1}{r}}\xi^{-1}}(x)I_{a,\xi}(x^{1/r},q^{-1}).$$
Then the functions $\tilde{I}_{a,\xi}(x,q^{-1})$ form a basis of solutions of the $q$-difference equation 
\begin{equation}
\label{intro:qdiffeqinfty}
\left[\prod_{k=1}^r(1-q^{rx\partial_x-k}) +(-1)^{r+N}xq^{\frac{r(r-1)}{2}+(r^{2}-N-1)x\partial_x} (1-q^{x\partial_x})^{N+1}\right] I =0,
\end{equation}
which is \eqref{intro:qdiffequ0} after the change of variable $x=Q^{-1}$.

\paragraph{Outline of the paper.}
In the first section, we recall the definition of the moduli space of $r$-spin curves, and we define the permutation-equivariant K-theoretic FJRW invariants.
In the second section, we define and compute the fake invariants by using Chiodo--Zvonkine's theorem \cite{chiodoTwistedRspinPotential2009}.
In the third section, we define and compute the spine CohFT, which is another building block of the FJRW invariants, related to $r$-spin curves with symmetries.
In the fourth section we use the Lefschetz formula to prove the adelic characterization theorem, which recursively determines the FJRW invariants.
In the last section we give an alternative description of the $J$-function in terms of untwisted invariants, and we use it to find a point $I_{\FJRW}$ in the image of the $J$-function.
\subsection*{Acknowledgment}
I am grateful to Y.P. Lee and Melissa Liu, whose suggestions lead to a more general version of the main result.
I would also like to thank Charles Doran, Mark Shoemaker, Yongbin Ruan, Yaoxiong Wen for their interest in this work, and  Xiaohan Yan for numerous fruitful discussions.
These results were obtained during my Ph.D. thesis at IMJ-PRG, and I am grateful to my advisor Alessandro Chiodo for his help all along the realization of this work.

\tableofcontents
\section*{Notations and conventions}
All schemes and stacks are of finite type over $\C$.
The Chow rings are taken with rational or complex coefficients
\begin{itemize}
\item $\overline{\cM}^r_{0,n}$, $\widetilde{\cM}^r_{0,n}$ : moduli space of $r$-spin curves, and moduli space of $r$-spin curves with trivialized marked points,
\item $\overline{\cM}^{r,s}_{0,n}$ : moduli space of $r$th roots of $\omega_{\log}^{\otimes s}$,
\item $\cL$ : universal $r$th root of $\omega_{\log}$,
\item $\cL_i$ : tautological line bundles over $\widetilde{\cM}^r_{0,n}$ ,
\item $\cI X$ : inertia stack of $X$,
\item $\ua,\ub,\ug$ : multi-indices,
\item $BG$ : classifying stack of the group $G$,
\item $\fe$ ramification index and cardinality of isotropy group,
\item $f_{(\xi)}$ : if $f$ is a Laurent polynomial in $q^{1/r}$, $f_{(\xi)}$ is the expansion of $f$ at $q^{1/r}=\xi$.
\end{itemize}
\section{Defining the invariants}

In this section we give the definition of the FJRW invariants in the permutation-equivariant case.
We first define a $K$-theoretic class $\Lambda$ over the moduli space of $r$-spin curves.
The symmetric group $S_n$ acts by permuting the marked points, and the class $\Lambda$ is $S_n$-equivariant in a natural way.
Thus, the cohomology groups $H^*(\Lambda)$ form $S_n$-modules.
To encode this representation in a way independent of $n$, we follow Getzler--Kapranov \cite{getzlerModularOperads1998} and Givental \cite{giventalPermutationEquivariantQuantumKtheory} and use Schur--Weyl duality. 
The representation $H^*(\Lambda)$ of $S_n$ is encoded in a symmetric function in infinitely many variables, which we take as a definition of the FJRW invariant.
More generally, we give a definition of the FJRW invariants with value in a chosen $\lambda$-ring $R$.
This procedure allows us to define generating functions for the FJRW invariants.
 
\subsection{Spin curves and the fundamental class}

\begin{defn}
Let $X$ be a Deligne--Mumford stack, and let $\cL$ be a line bundle over $X$.
We say that $\cL$ is \emph{representable} if, for any geometric point $x$ of $X$, the induced representation $\Aut(x) \to \G_m$ is faithful.

If $X$ is of finite type, $\cL$ is representable if and only if the induced morphism $X\to B\G_m$ is representable (see \cite[4.4.3.]{abramovichCompactifyingSpaceStable2001}).
\end{defn}

\begin{defn}
Let $n\geq 3$ and $r\geq 1$ be integers.
An $r$-spin curve with $n$ marked points is the data of $(C,(\Sigma,\hdots,\Sigma_n),\cL,\alpha)$, where 
\begin{itemize}
\item $(C,(\Sigma_1,\hdots, \Sigma_n))$ is a stable balanced twisted curve in the sense of \cite{abramovichTwistedBundlesAdmissible2003},
\item $\cL$ is a representable line bundle over $C$,
\item $\alpha: \cL^{\otimes r}\to \omega_{\log}$ is an isomorphism.
\end{itemize}

The moduli stack of $r$-spin curves $\overline{\cM}^r_{g,n}$ classifies families of $r$-spin curves.
There is a universal curve $\pi : \overline{\cC}^r_{g,n}\to \overline{\cM}^r_{g,n}$, and a universal line bundle $\cL_{g,n} \to \overline{\cC}^r_{g,n}$.
The marked points $\Sigma_i \subset \overline{\cC}^r_{g,n}$ are closed substacks of the universal curve.

We also define the stack $\widetilde{\cM}^r_{0,n}$ classifying $r$-spin curves together with a section $\sigma_i$ at each marked point, namely,
$$\widetilde{\cM}^r_{0,n} = \Sigma_1\times_{\overline{\cM}^r_{g,n}} \hdots \times_{\overline{\cM}^r_{g,n}} \Sigma_n.$$

The forgetful map $p:\widetilde{\cM}_{g,n}^r\to \overline{\cM}_{g,n}^r$ is a $\bmu_r^n$-gerbe.
The universal curve of $\widetilde{\cM}^r_{0,n}$ is $\widetilde{\cC}_{g,n}^r = \overline{\cC}^r_{g,n}\times_{\overline{\cM}^r_{g,n}} \widetilde{\cM}^r_{g,n}$, and we still denote by $\pi : \widetilde{\cC}_{g,n}^r\to \widetilde{\cM}_{g,n}^r$ the projection.
The stack $\widetilde{\cM}^r_{0,n}$ is equipped with with the tautological line bundles $\cL_i = \sigma_i^* \omega_\pi$, where $\omega_\pi$ is the dualizing line bundle.
\end{defn}
\subsubsection{Multiplicities}
At each marked point $\Sigma_i$ of an $r$-spin curve $C$, the stabilizer is canonically isomorphic to $\mu_{\fe_i}$, for some $\fe_i|r$.
The line bundle $\cL_{|\Sigma_i}$ is a representation of $\mu_{\fe_i}$, given by some element $d_i\in \Z/\Z_{\fe_i}$, with $d_i\wedge \fe_i=1$.
We refer to $\frac{d_{i}}{\fe_i}\in \Q/\Z$ as the \emph{multiplicity} of $\cL$ at $\Sigma_i$, and we denote it by $\mult_{\Sigma_i}(\cL)$.
By construction, $\mult_{\Sigma_{i}}(\cL)$ is actually an element of $\frac{1}{r}\Z/\Z$.
For a multi-index $\ua \in \left(\Z/r\Z\right)^{n}$, we denote by $\overline{\cM}^{r}_{0,\ua}$ the open and closed substack of $\overline{\cM}^{r}_{0,n}$ parametrizing $r$-spin curves with multiplicity $\frac{\ua}{r}$.
The marked points with multiplicity $0$ are called \emph{broad marked points}.

\subsubsection{Evaluation maps}
Recall \cite{abramovichGromovWittenTheoryDeligneMumford2008} that for a given stack $\cX$, $\cI_{\bmu_{k}}\cX$ classifies representable maps from a \emph{trivialised} $\bmu_{k}$-gerbe $X\times B\bmu_{k}$ to $X$.
The \emph{cyclotomic inertia stack} is $\cI_{\bmu}\cX = \bigsqcup_{k\geq  1}\cI_{\bmu_{k}}\cX$. 
\begin{rem}
	Over $\C$, there is a canonical isomorphism $\cI_{\bmu}\cX \simeq \cI \cX$ following from the isomorphism $\bmu_{k}\simeq \Z/k\Z$.
	\end{rem}

By construction, the gerbes $\Sigma_{i}$ are canonically trivialised over $\widetilde{\cM}^r_{\ua}$, and the restrictions $(\omega_{\log})_{|\Sigma_{i}}$ are canonically trivial by the residue map.
Thus, $\cL_{|\Sigma_{i}}$ is an $r$th root of the trivial bundle, and defines a morphism $\Sigma_{i}\to B\bmu_{r}$.
\begin{defn}
	The \emph{evaluation maps} are the morphisms $$\ev_{i} : \widetilde{\cM}^{r}_{0,n} \to \cI B\mu_{r}$$
	associated to the trivial gerbe $\Sigma_{i}$ and the morphism $\Sigma_{i} \to B\mu_{r}$ as discussed above.
\end{defn}
The automorphism $\xi \mapsto \xi^{-1}$ of $\mu_{r}$ induces an automorphism $\iota$ of $\cI B\bmu_{r}$.
In order to glue $r$-spin curves, we need the twisted evaluation maps.
\begin{defn}
The twisted evaluation maps are 
$$\ev_{i}^{\vee} = \iota \circ \ev_{i} .$$
\end{defn}

\begin{rem}
	We can give an  explicit description of $\ev_{i}$ and $\ev_{i}^{\vee}$ as follows.
	The pullback $\sigma_{i}^{*}\cL$ is canonically an $r$th root of the trivial line bundle, and defines a map $\widetilde{\cM}^r_{\ua} \to B\bmu_{r}$.
On the other hand, $\cI B\bmu_{r}$ is a disjoint union of $r$ copies of $B\bmu_{r}$ $$\cI B\bmu_{r} \simeq \bigsqcup_{a=0}^{r-1}B\mu_{r}.$$
	The $i$th evaluation map $\ev_{i}$ sends $\widetilde{\cM}^r_{\ua}$ to the $a_{i}$th copy of $B\bmu_{r}$ via the map above, while $\ev_{i}^{\vee}$ goes to the $-a_{i}$th copy via the map induced by $\cL^{\vee}$.
\end{rem}

\subsubsection{Fundamental class}
In order to define the fundamental class, we need $R^1\pi_*\cL$ to be a vector bundle.
This is the case, for example, if $H^{0}(C,\cL)=0$ for every $r$-spin curve $(C,\cL)$ corresponding to a closed point of $\overline{\cM}_{0,\ua}^{r}$.
However, when a spin curve has at least 2 broad points, it may happen that $H^{0}(C,\cL)$ is non-zero.
To fix this, we twist $\cL$ by the divisor of broad marked points before pushing forward to $\overline{\cM}^r_{0,n}$.  
\begin{lem}
Let $E\subset \overline{\cC}^r_{\ua}$ be the divisor of broad marked points, that is, $E = \bigsqcup_{a_i=0} \Sigma_i$.
Then we have $\pi_*\cL_{\ua}(-E)=0$, and $R^1\pi_*\cL(-E)$ is a vector bundle over $\overline{\cM}^r_{0,\ua}$.
\end{lem}

\begin{proof}
See \cite[lem 4.1.1.]{chiodoLandauGinzburgCalabiYauCorrespondence2010}.
\end{proof}

\begin{rem}
This twist has a mild impact on the Chern classes.
Indeed, we have the exact sequence $$0\to \cL(-E) \to \cL \to \cL_{|E} \to 0.$$ 
For a broad marked point $\Sigma_i$, we have $p_*\left( \cL_{|\Sigma_i}\right) = \sigma_i^* \cL$, which shows that $c\left( p_*\left( \cL_{|\Sigma_i}\right)\right) = 1$ (in $\mathrm{A}^\bullet_\Q$).
Thus we have $c(R\pi_*\cL)=c(R\pi_*\cL(-E))$, and $\ch\left(R\pi_*\cL\right)=\ch\left(R\pi_*\cL(-E)\right) +m$, where $m$ is the number of broad marked points.
\end{rem}

\begin{prop}
Let $S_n$ act on $\widetilde{\cM}^r_{0,n}$ and $\overline{\cM}^{r}_{0,n}$ by permuting of the marked points.
Then the sheaf $R^1\pi_*\cL\in \mathrm{Coh}\left(\overline{\cM}^r_{0,n}\right)$ is naturally $S_n$-equivariant, and so is its pullback to $\widetilde{\cM}^r_{0,n}$.
\end{prop}

\begin{rem}
Note that the different connected components of $\widetilde{\cM}^r_{0,n}$ may be permuted by the $S_n$-action.
\end{rem}

\begin{defn}[The fundamental class]
The fundamental class of the Fermat polynomial $W=\sum_{i=0}^{N}X_{i}^{r}$ is the $S_n$-equivariant K-theoretic class  $$\Lambda_n = \left( \lambda_{-1}\left(R^1\pi_*\cL(-E)\right)\right)^{\otimes N+1}\in K^0_{S_n}(\overline{\cM}^r_{0,n}). $$
More generally, we define
$$\Lambda_n(s) = \left( \lambda_{-s}\left(R^1\pi_*\cL(-E)\right)\right)^{\otimes N+1}\in K^0_{S_n}(\overline{\cM}^r_{0,n})\llbracket s \rrbracket.$$ 
\end{defn}

\subsection{A symplectic space}

Following Givental (\cite{giventalSymplecticGeometryFrobenius2004,giventalHirzebruchRiemannRoch2011}) we define an infinite-dimensional symplectic space $\cK$, referred to as the loop space.
This space is equipped with a natural polarization $\cK = \cK_+\oplus \cK_-$. 

\begin{defn}
The \emph{state space} of the FJRW theory is $K^0\left( \cI B\bmu_r\right)_\C$.
Since $\cI B\bmu_r \simeq \bigsqcup_{a=0}^{r-1} B\bmu_r$, there is an isomorphism of \emph{vector spaces}
$$K^0\left( \cI B\bmu_r\right)_\C \simeq \C\left[ \Z_r\right] \otimes_\C \C\left[\widehat{\bmu}_r \right].$$

There are two natural basis for this vector space.
For $a,l\in \Z/r\Z$, let $\phi_{a,l}$ be the character $\zeta \mapsto \zeta^{l}$ on the $a$th copy of $B\bmu_{r}$, and for $\xi\in\bmu_{r}$, we define $$\phi_{a,\xi} = \frac{1}{r} \sum_{l\in\Z_{r}} \xi^{-l}\phi_{a,l}.$$
Then the sets $\{\phi_{a,l}\}_{a,l\in\Z_{r}}$ and $\{\phi_{a,\xi}\}_{\substack{a\in\Z_{r}\\\xi\in \bmu_{r}}}$ are both basis of $K^{0}(\cI B\bmu_{r})$.

The state space is equipped with the orbifold pairing twisted by the fundamental class
\begin{equation}
	\scal{\phi_{a,\xi}, \phi_{a',\xi'}} = \left\{ \begin{matrix}\frac{1}{r}\delta_{a,-a'}\delta_{\xi,\xi'} & \textrm{ if } a \neq 0, \\ \frac{(1-\xi s)^{N+1}}{r} \delta_{0,a'} \delta_{\xi,\xi'} & \textrm{ otherwise.} \end{matrix} \right.
\end{equation}
The dual of the element $\phi_{a,l}$ (resp. $\phi_{a,\xi}$) with respect to this pairing is denoted by $\phi^{a,l}$ (resp. $\phi^{a,\xi}$).
Finally, the Adams operation of K-theory are the ring morphisms
$$\Psi^m(\phi_{a,\xi})= \sum_{\zeta^m=\xi}\phi_{a,\zeta}.$$
\end{defn}

For each element $a \in \Z/r\Z$, let $\fe(a)$ be the order of the subgroup generated by $a$ in $\Z/r\Z$.
If $\cL$ has multiplicity $\frac{a}{r}$ at a marked point $x_i$, then its automorphism group $\Aut(x_i)$ is isomorphic to $\Z/\fe(a)\Z$. 
We decompose $\C[\Z/r\Z]$ according to this order $\fe$ :
\begin{equation}
\C[\Z/r\Z] = \bigoplus_{\fe|r} V_\fe,
\end{equation}
with $V_\fe = \bigoplus_{\fe(a)=\fe} \C\cdot \phi_a$.

\begin{defn}
\label{defn:loopspace}
The loop space of the FJRW theory is the space of rational functions
\begin{equation}
\cK = \bigoplus_{\fe|r} V_\fe(q^{\frac{1}{\fe}})\otimes \C[\widehat{\bmu}_r].
\end{equation}

The loop space $\cK$ is equipped with the symplectic form $\Omega$, whose restriction to $V_\fe(q^{\frac{1}{\fe}})\otimes \C[\widehat{\bmu}_r] $ is 
\begin{equation}
\Omega(f,g)=\left[\Res_{q^{1/\fe}=0} + \Res_{q^{1/\fe}=\infty}\right] \scal{f(q^{1/\fe}),g(q^{-1/\fe})}\frac{dq^{1/\fe}}{q^{1/\fe}}.
\end{equation}

We define a polarization of this symplectic vector space by setting 
\begin{align*}
\cK_+ &= \bigoplus_\fe V_\fe\left[q^{\frac{1}{\fe}},q^{\frac{-1}{\fe}}\right] &
\cK_- &= \left\{ f\in \cK| f(0)\neq \infty \textrm{ and } f(\infty) =0\right\}.
\end{align*}
\end{defn}
\begin{notation}
	The unit of $K^{0}(\cI B\bmu_{r})$ for the orbifold tensor product is $\phi_{1,0}$, which we will omit in elements of $\cK$.
	Thus, elements of the form$f(q^{1/r})\phi_{1,0}$ will be simply denoted by $f$.
\end{notation}

\subsection{The invariants}
\begin{defn}
\label{defn:S_nmodules}
Let $t$ be an element of $\cK_+$, and $i\in \{1,\hdots, n\}$
We introduce the class $t(\cL_i)\in K^0\left(\widetilde{\cM}^r_{0,n}\right)$, which is defined on elementary tensors by $$(\phi_{a,\xi}) q^{j/\fe(a)} \mapsto \ev_i^*(\phi_{a,\xi}) \otimes \cL_i^{\otimes j} .$$
For any $t\in \cK_+ $, the class $\Lambda_{n}(s)\otimes \bigotimes_{i=1}^nt(\cL_i)$ is naturally  an $S_n$-equivariant class, and its cohomology groups form an $S_n$-module denoted by $$\left[ t(\cL_1),\hdots, t(\cL_n)\right]_{n} = H^*\left( \widetilde{\cM}^r_{0,n} ; p^*\Lambda_n(s) \otimes \bigotimes_{i=1}^nt(\cL_i)\right).$$

More generally, let $n=k_1+\hdots +k_s$ be a partition of $n$, and let $H\subset S_n$ be the subgroup $S_{k_1}\times \hdots \times S_{k_s}$.
Let us denote $x_{i,k}$ ($i\in \{1,\hdots, s\}$ and $k\in \{1,\hdots, k_s\}$) the $n$ marked points.
Then, for a sequence of inputs $t^{(1)},\hdots, t^{(s)}\in \cK_+$, the cohomology groups $$H^*\left( \widetilde{\cM}^r_{0,n} ; p^*\Lambda_n(s) \otimes \bigotimes_{k=1}^s \bigotimes_{l=1}^{k_s}t^{(k)}(\cL_{k,i})\right)$$ are $H$-modules, denoted by $$\left[ t^{(1)}(\cL_{1,1}),\hdots, t^{(1)}(\cL_{1,k_1});\hdots,t^{(s)}(\cL_{s,k_s})\right]_{n}.$$
\end{defn}

In order to define a generating function for these $S_n$-modules, we use the ring of symmetric function (\cite{getzlerModularOperads1998}, \cite{giventalPermutationEquivariantQuantumKtheory}).
This  allows us to encode representations of $S_n$ for various $n$ in a single ring.

More generally, let $R$ be a $\lambda$-ring over $\C$.
We assume that $R$ is equipped with the $\scrI$-adic topology for an ideal $\scrI$ of $R$, such that 
\begin{itemize}
\item $R$ is Hausdorff,
\item for all $m\geq 0$, $\Psi^m(\scrI) \subset \scrI^m$.
\end{itemize} 
In that case, the completion $\widehat{R}$ of $R$ is also a $\lambda$-ring.
The main examples are the ring of symmetric functions, and $\C[X]$, and $\scrI$ is the ideal of functions with constant term equal to 0.

We extend the scalar to $R$ in \Cref{defn:loopspace}, and complete the resulting ring with respect to the $\scrI$-adic topology (see \cite[Appendix B]{coatesComputingGenusZeroTwisted2009} for a detailed construction of the loop space).
In particular, $\cK_+$ is made of functions $t$ which, modulo any power of $\scrI$, are Laurent polynomials.

We now come to the definition of the FJRW invariants.
\begin{defn}
\label{defn:invariants}
We keep the notations of \Cref{defn:S_nmodules}.
For any elements $\nu_1,\hdots, \nu_s \in R$, we define

\begin{equation}
\label{defn:equation:invariantpermeq}
\scal{t^{(1)}(\cL_{1,1})\otimes \nu_1, \hdots ,t^{(s)}(\cL_{s,k_s})\otimes \nu_s}_{0,n}^{H}  =  \frac{1}{\prod_i k_i!}\sum_{h\in H}\tr_h \left[ t^{(1)}(\cL_{1,1}),\hdots,t^{(s)}(\cL_{s,k_s})\right]_{n} \prod_{i=1}^s \prod_{j=1}^\infty \Psi^r(\nu_i)^{l_j(h)}
\end{equation}

where $l_j(h)$ is the number of cycles of length $j$ in $h$.
\end{defn}

\begin{rem}
In the case where $H=S_n$, and $R$ is the ring of symmetric functions, \eqref{defn:equation:invariantpermeq} yields the symmetric function associated to the $S_n$-module $\left[t(\cL_1),\hdots,t(\cL_n)\right]$.
When $R=\Q[x]$, the FJRW invariants $$\scal{t(\cL_{1})\otimes x, \hdots ,t(\cL_{n})\otimes x}_{0,n}^{S_{n}}$$ correspond to the (virtual) dimension of the fixed subspace.
\end{rem}

\begin{rem}[Vanishing]
\label{rem:vanishing}
With the same notations as in the previous definition, let us choose $t^{(i)} = q^{j/\fe(a)}\phi_{a,l}$.
Then, the invariant vanishes unless $l a=j\mod \fe_i$.

Indeed, consider the forgetful map $p:\widetilde{\cM}^r_{0,n} \to \overline{\cM}^r_{0,n}$.
The line bundle $\cL_i^j\otimes\ev_i^*[d]$ carries the $\bmu_{\fe_i}$ representation $\zeta\mapsto \zeta^{-j}\zeta^{ad}$, where $\bmu_{\fe_i}$ is the group of $2$-automorphisms of the section $\sigma_i$.
Then, $p_*(\cL_i^j\ev_i^*[d])=0$ if this representation is non-trivial.
\end{rem}

\begin{prop}[poly-linearity]
For $t,t'\in \cK_+$, and $\nu,\nu'\in R$, we have
\begin{equation}
\scal{t\otimes \nu + t'\otimes \nu',\hdots, t\otimes \nu + t'\otimes \nu'}^{S_n}_{0,n} = \sum_{k+l=n}\scal{t\otimes \nu,\hdots, t\otimes \nu ,t'\otimes \nu',\hdots, t'\otimes \nu'}^{S_k\times S_l}_{0,n}.   
\end{equation}
We use this formula to extend \Cref{defn:invariants} to inputs $t\in \cK_+$.
\end{prop}
\begin{proof}
See \cite{giventalPermutationEquivariantQuantumKtheory} example 5.
\end{proof}

\begin{hypo}
From now on we will assume that the input $t$ belongs to $\scrI \cK_+$.
This ensures that the following formal series are well-defined.
\end{hypo}

\begin{defn}
The genus-0 permutation-equivariant potential is the formal function $\cF_0$, defined over $\cI\cK_+$ by
\begin{equation}
\cF_0(t)=\sum_{n\geq 0} \scal{t(\cL_1),\hdots, t(\cL_n)}^{S_n}_{0,n}
\end{equation}
We also consider the mixed potential
\begin{equation}
\cF_0(x,t)=\sum_{n+k\geq 3} \frac{1}{k!}\scal{x(\cL_1),\hdots, x(\cL_k),t(\cL_{k+1}),\hdots, t(\cL_{n+k})}^{S_n}_{0,k+n}.
\end{equation}
\end{defn}

We now introduce the $J$-function which, up to a translation, is the differential of the mixed potential with respect to the first variable.

\begin{defn}
\label{defn:J}
The $J$-function is the formal function $\scrI\cK_+ \to \cK$ defined by
\begin{equation}
J(t) = 1-q + t + \sum_{\substack{n\geq 2\\a\in \Z_r\\ \xi \in \mu_r}} \phi^{a,\xi} \scal{\frac{\phi_{a,\xi}}{1-q^{\frac{1}{\fe(a)}}\cL_0},t(\cL_1),\hdots, t(\cL_{n})}^{S_n}_{0,n+1}.
\end{equation}
\end{defn}

\section{The fake theories}

This section is devoted to the definition and computation of the so-called fake theories, which can be seen as building block for the K-theoretic invariants.
Using the theory of twisted invariants developed by Coates, Givental, and Tonita in \cite{giventalSymplecticGeometryFrobenius2004,coatesQuantumRiemannRoch2007,tonitaTwistedOrbifoldGromovWitten2014}, and a theorem by Chiodo and Zvonkine \cite{chiodoTwistedRspinPotential2009}, we are able to fully compute these fake theories.
\subsection{The fake invariants}
\label{subsection:fakeinvariants}
\begin{defn}
	\label{defn:fakeinvariants}
Let $\cA_\xi,\cB$, and $\cC $ be invertible multiplicative classes, and let $Z$ denote the singular locus in the universal curve $\overline{\cC}_{0,n}^r \to \overline{\cM}^r_{0,n}$.
Define the following classes 
\begin{align*}
\cA_{\xi,n}(\ua) &= r\cA_\xi(R\pi_*\cL_{\ua}(-E)) \in H^*\left(\overline{\cM}^r_{\ua}\right)\\
\cB_{0,n} &= \cB\left(\pi_*\left(\omega_{\log}^{-1}-1\right) \right)\\
\cC_{0,n} &= \cC\left(\pi_*\cO_{Z}\right).
\end{align*}
The fake invariants are defined by
\begin{equation}
\scal{\phi_{a_1,\xi_1}\cL_1^{k_1},\hdots, \phi_{a_n,\xi_n}\cL_n^{k_n}}^\fake_{0,n} = \left\{ \begin{matrix}\int_{\overline{\cM}^r_{0,\ua}}\cA_{\xi,n}(\ua)\cB_{0,n}\cC_{0,n} \prod_{i=1}^n \ch(\cL_i^{k_i}) & \textrm{ if } \xi_i = \xi \forall i,\\ 0 & \textrm{ otherwise.} \end{matrix} \right.
\end{equation}

\end{defn}

\begin{rem}
Notice a slight abuse of notation in the definition above.
Indeed, the tautological line bundles $\cL_i$ do not live on $\overline{\cM}^r_{0,n}$, but rather on $\widetilde{\cM}^r_{0,n}$.
Thus, in the definition above, $\ch(\cL_i)$ should be interpreted as $e^{\frac{\psi_i}{\fe_i}}$, where $\psi_i$ is the usual $\psi$-class pulled back from $\overline{\cM}_{0,n}$.
\end{rem}

For the rest of this article, we choose 
\begin{align*}
\cB(L) &= \td^{-1}(L)\\
\cC(L)&= \td^{-1}(L^{\vee})
\end{align*}
With this definition, the fake invariants become 
$$\scal{\phi_{a_1,\xi}\cL_1^{k_1},\hdots, \phi_{a_n,\xi}\cL_n^{k_n}}^\fake_{0,n} = \int_{\overline{\cM}^r_{0,\ua}}\cA_{\xi,n}(\ua) \prod_{i=1}^n \ch(\cL_i^{k_i}) \td(\cT),$$
where $\cT$ is the tangent space.

\subsection{Fake $J$-functions}
Following the work of Givental \cite{giventalSymplecticGeometryFrobenius2004} and Tonita \cite{tonitaTwistedOrbifoldGromovWitten2014} we organize these invariants in the so-called $J$-function.

\begin{defn}
\label{defn:kfake}
Let $\cK^\fake$ be the vector space
\begin{equation}
\cK^{\fake} = \bigoplus_\fe V_\fe\otimes \C[\widehat{\bmu}_r]\otimes \C\llbracket q^{1/\fe}-1,(q^{1/\fe}-1)], 
\end{equation}
equipped with the symplectic form
\begin{equation}
\Omega^{\fake}(f,g) = r\Res_{q^{1/r}=1} \scal{f(q^{-1/r}),g(q^{1/r})}^\cA\frac{dq^{1/r}}{q^{1/r}},
\end{equation}
where the inner product $\scal{..}^{\cA}$ is defined over $\C[\Z_r]\otimes \C[\widehat{\mu}_r]$ by
\begin{align*}
	\scal{\phi_{a,\xi} , \phi_{b,\zeta}}^\cA &= \left\{ \begin{matrix}\delta_{a,-b}\delta_{\xi,\zeta} & \textrm{ if } a\neq 0,\\ \delta_{0,b} \delta_{\xi,\zeta} \cA_\xi^{-1}(\cO) & \textrm{otherwise.}\end{matrix} \right.
\end{align*}
The dual basis is denoted by $\{\phi^{a,\xi}\}$.
\begin{rem}
	There is a slight conflict of notation with the previous section, because the inner products on $\cK$ and $\cK^{\fake}$ differ by a factor $r$.
	Thus, the dual basis $\phi^{a,\xi}$ also differs by a factor $r$ depending on wether we consider them as elments in $\cK$ or $\cK^{\fake}$.
	This difference is compensated in the $J$-function by the fact that the fake invariants also have a factor $r$ (\Cref{defn:fakeinvariants}) compared to the K-theoretic $J$-function.
\end{rem}
We equip this symplectic space with the polarization
\begin{align*}
\cK^{\fake}_+ &= \bigoplus_\fe V_\fe\otimes \C[\widehat{\bmu}_r]\llbracket q^{1/\fe} -1\rrbracket,\\
\cK^{\fake}_- &= \bigoplus_\fe V_\fe \otimes \C[\widehat{\bmu}_r] [ (q^{1/\fe} -1)^{-1}].
\end{align*}
The potential of the fake theory is the formal function defined on $\cK^{\fake}_+$ by 
\begin{equation}
\cF^{\fake}(t)= \sum_{n\geq 3} \frac{1}{n!}\scal{t(\cL_1),\hdots, t(\cL_n)}^\fake_{0,n}.
\end{equation}
The fake $J$-function is the shifted graph of the differential of $\cF$ inside $\cK^\fake$
\begin{equation}
	J(t)=1-q +t + \sum_{\substack{n\geq 2\\a\in \Z_r}} \frac{ \phi^{a,\xi}}{\fe(a)n!} \scal{\frac{\phi_{a,\xi}}{1-q^{1/\fe(a)}\cL_0},t(\cL_1),\hdots,t(\cL_n)}^\fake_{0,n+1}.
\end{equation}
\end{defn}

\subsection{Lagrangian cones}
The image of the $J$-function is a Lagrangian cone in $\cK^\fake$, which can be explicitly computed, as we now explain.
The collection of classes $\cA_{\xi,n}(\ua)$ form the genus-0 part of a CohFT over the state space $\C[\Z_r]\otimes \C[\widehat{\bmu}_r]$.
Its associated Lagrangian cone $L^\cA$ (see \cite{giventalSymplecticGeometryFrobenius2004}) lies in the symplectic space $(\cH^\cA,\Omega^\cA$) given by 
\begin{align*}
\cH^{\cA} &= V\llbracket z, z^{-1}] & \Omega^{\cA}(f(z),g(z)) &= \Res_{z=0} \scal{f(-z),g(z)}^{\cA}dz
\end{align*}
For $\cA_\xi=1$, the resulting cone is called the \emph{untwisted cone} $L^\un$, and can be easily deduced from the cohomological $J$-function of a point.

\begin{prop}[Chiodo--Zvonkine \cite{chiodoTwistedRspinPotential2009}]
\label{prop:Delta}
Let $w_\xi(z) = \sum_{g\geq 0}w_{d,\xi}z^d\in \C\llbracket z \rrbracket$ be  power series, and let $\cA_\xi$ be the multiplicative classes
\begin{equation}
\cA_\xi(E)= \exp\left( \sum_{d\geq  0} w_{d,\xi} \ch_d(E)\right).
\end{equation}
Let $\Delta$ be the operator acting on $\cH^{\cA}$ such that for all $0\leq a \leq r-1$ we have 
\begin{equation}
	\Delta \left( \phi_{a+1,\xi} \right)= \exp\left( \sum_{d} w_{d,\xi} \frac{B_{d+1}(\frac{a+1}{r})}{(d+1)!}z^d\right)\phi_{a+1,\xi}.
\end{equation}
Then we have 
\begin{equation}
\Delta L^\un = L^\cA.
\end{equation}
\end{prop}
\begin{rem}
The shift in the definition of $\Delta$ happens because of the twist by the divisor of broad marked points $\cL(-E)$.
\end{rem}
Finally, it is a consequence of \cite{tonitaTwistedOrbifoldGromovWitten2014} that the cones $L^\fake$ and $L^A$ coincide.

\begin{prop}[\cite{tonitaTwistedOrbifoldGromovWitten2014}]
Let $\ch$ be the morphism
\begin{align*}
\ch : \cK^\fake &\to \cH^\cA \\
q^{j/r}\phi_{a,\xi} &\mapsto e^{jz/r}\phi_{a,\xi}.
\end{align*}
Then we have 
\begin{equation}
L^\fake = \ch^{-1}\left(L^\cA\right).
\end{equation}
\end{prop}
Since $\ch$ is an isomorphism, we identify $L^\fake$ and $L^\cA$, and write $L^\fake=L^\cA$.
We now apply these results to the classes $$\cA_\xi(E) = \ch\left( \lambda_{-s\xi}E \right)^{-(N+1)}.$$

\begin{prop}
For $\cA_\xi$ as above, we have $$w_{d,\xi} = (N+1)\sum_{k\geq 1} \frac{s^k\xi^{k}k^d}{k},$$
and $$\Delta (\phi_{a,\xi}) = \left\{ \begin{matrix} \exp\left( (N+1)\sum_{k\geq 1} \frac{s^k \xi^{k}}{k} \frac{q^{\frac{ka}{r}}}{q^{k}-1} \right) \phi_{a,\xi}& \textrm{ if } a\neq 0 ,\\
\exp\left((N+1) \sum_{k\geq 1} \frac{s^k \xi^{k}}{k} \frac{q^{k}}{q^{k}-1} \right)\phi_{0,\xi}  & \textrm{ otherwise.}  \end{matrix}\right.$$
\end{prop}
 
\subsection{An extension of the fake theory} 
\label{defn:krm}
In order to deal with permutations of the marked points in the next section, we need to generalize slightly the previous definition to include $rm$th roots of $\omega_{\log}^{\otimes m}$, and $rm$th roots of unity.
Indeed, an $r$-spin curve with an automorphism of order $m$ naturally yields an $rm$-root of $\omega_{\log}$ on the quotient curve (see \Cref{subsection:Spine contribution}).
The moduli space of $rm$th roots of $\omega_{\log}^{\otimes m}$ is $\overline{\cM}_{0,n}^{rm,m}$, and has a forgetful map $\epsilon : \overline{\cM}^{rm,m}_{0,n} \to \overline{\cM}_{0,n}$.
The universal curve $\overline{\cC}^{rm,m}_{0,n}$ carries the universal $rm$th root $\cL$.
The multiplicity of $\cL$ at a marked point is now an element of $\left(\frac{1}{rm}\Z\right)/\Z$.

\begin{defn}
For multiplicative classes $\cA_\xi = \exp\left( \sum w_{d,\xi}\ch_d\right)$, $\xi\in \bmu_{rm}$, and $\ua$ a multi-index, we define 
\begin{equation}
\cA_{\xi,n}(\ua)= rm \epsilon_* \cA_\xi\left( R\pi_*\cL_{\ua}(-E)\right),
\end{equation}
where $\cL_{\ua}$ is the universal $rm$th root of $\omega_{\log}^m$, with multiplicity $\ua$.
These classes form a genus-0 CohFT over the state space $W=\C[\Z/{rm}\Z]\otimes \C[\widehat{\bmu}_{rm}]$.
We write $W$ as the direct sum $$\C[\Z/rm\Z] = \bigoplus_{\fe|rm} W_\fe,$$
where $W_\fe$ is spanned by the basis elements $\phi_a$ such that the order of $a$ in $\Z/rm\Z$ is $\fe$.
\end{defn}

For the remaining part of this article, we fix $$\cA_\xi(E) = \ch\left(\lambda_{-s\xi}E\right)^{-(N+1)},$$ for $\xi \in \bmu_{rm}$. 
By \cite{chiodoTwistedGromovWittenRspin2007}, the associated Lagrangian cone of this CohFT is equal to $ \Delta L^\un$, with 
\begin{align*}
\Delta (\phi_{a+1,\xi})&= \exp\left( (N+1) \sum_{k\geq 1} \frac{s^k\xi^{k} k^d}{k(d+1)!}z^d B_{d+1}\left(\frac{a+1}{rm}\right)\right) \phi_{a+1,\xi}\\
&=\exp\left( (N+1)\sum_{k\geq 1} \frac{(s\xi)^{k}}{k} \frac{q^{\frac{k(a+1)}{rm}}}{q^{k}-1}\right) \phi_{a+1,\xi}  \textrm{ for } 0\leq a \leq rm-1
\end{align*}
We extend the fake invariants to $\C[\Z_{rm}] \otimes \C[\widehat{\bmu}_{rm}]$ by setting 
\begin{equation}
\scal{\phi_{a_1}\otimes e_{\xi_1}\cL_1^{k_1},\hdots, \phi_{a_n}\otimes e_{\xi_n}\cL_n^{k_n}}^\fake_{0,n} = \left\{ \begin{matrix}\int_{\overline{\cM}^r_{0,\ua}}\cA_{\xi,n}(\ua)\cB_{n}\cC_{n} \prod_{i=1}^n \ch(\cL_i^{k_i}) & \textrm{ if } \xi_i = \xi \forall i,\\ 0 & \textrm{ otherwise.} \end{matrix} \right.
\end{equation}
%
The associated Lagrangian cone $L^\fake$ lies in the symplectic space $$\cK_{rm} = \C\left[\Z_{rm}\right]\otimes \C\left[\widehat{\mu}_{rm}\right] \left\llbracket q-1,(q-1)^{-1}\right],$$
By \cite{tonitaTwistedOrbifoldGromovWitten2014}, the polarization of $\cK_{rm}$ is given by
\begin{align*}
\left( \cK_{rm}\right)_+ &= \bigoplus_\fe W_\fe\otimes \C\left[\widehat{\bmu}_r\right] \left\llbracket 1-q^{1/\fe} \right\rrbracket\\
\left( \cK_{rm}\right)_- &= \bigoplus_\fe W_\fe\otimes \C\left[\widehat{\bmu}_r\right] \left[ \left(1-q^{1/\fe}\right)^{-1} \right],
\end{align*}
With this choice of polarization, we have $$L^\fake = \Delta L^\un.$$
By a slight abuse of notation, we kept the notation $L^{\fake}$ for the Lagrangian cone of the extended fake theory.
This abuse of notation is justified by the following proposition.
\begin{prop}
Let $\Phi_0$ be the inclusion morphism
\begin{align*}
\Phi_0 : \cK^{\fake} &\to \cK_{rm}\\
\phi_{a,\xi} &\mapsto \phi_{ma,\xi}.
\end{align*}
Then $\Phi_0$ is an isomorphism of polarized symplectic spaces onto its image, and we have $\Phi_0(L^\fake) \subset L^\fake$.
\end{prop}

\section{The spine CohFT}
\label{section:spinecohft}
This section is devoted to the definition and computation of the spine CohFT.
This CohFT is designed to reproduce the moduli space of heads (see \ref{defn:heads}), and to recover the spine contribution of \Cref{adelic}.
Let us briefly sketch how the spine CohFT arises.
Given an $r$-spin curve $(C,\cL)$ with an automorphism $g\in \Aut(C)$ of order $m$ and an isomorphism $\phi:g^*\cL\to \cL$ (compatible with the spin structure) there is a line bundle $\bar{\cL}$ on the quotient curve $D=\left[C/g\right]$ constructed by descent.
The line bundle $\bar{\cL}$ is canonically an $rm$th root of $\omega_{\log}^{\otimes m}$ on $D$.
Thus, $D$ is equipped with this $rm$th root $\bar{\cL}$, and the $m$th root $T$ of the trivial bundle corresponding to the $\Z_m$-cover $C\to D$.
We take this situation as a definition, and we consider the moduli space $ \overline{\cM}_{0}^{rm}(B\Z_m;\ua,\ub)$ parametrizing curves with an $rm$th root $\overline{\cL}$ of $\omega_{\log}^{m}$ and an $m$th root $T$ of the trivial line bundle. 

Actually, only a small part of the spine CohFT will be relevant to our study.
Indeed, we only have to consider the case where the automorphism $g$ fixes two marked points and acts freely on the remaining marked points.
This is equivalent to asking that $T$ has trivial multiplicity at every marked point except two, where the multiplicity is co-prime to $m$.
The remaining part of the CohFT plays no role, so we can safely assume that each $b_{i}$ is either $0$, or prime to $m$.

\subsection{Stable maps to $B\Z_m$ and roots of the trivial bundle}
We recall the well-known correspondence between cyclic covers and roots of the trivial line bundle.
For a multi-index $\ug \in (\Z/m\Z)^n$, the space $\overline{\cM}_{g,n}(B\Z_m,\ug)$ parametrizes stable maps to $B\Z_m$ with holonomy $\ug$.
This stack admits an other description in terms of $m$th roots of the trivial bundle.
For a curve $D$, a stable map $D\to B\Z_{m}$ is given by a $\Z/m\Z$ cover $p:C\to D$, with holonomy $\ug$ at the marked points.
Let $\sigma$ be the canonical generator of $\Z_m$, and $\zeta=e^{\frac{2i\pi}{m}}$.
At a marked point $x\in p^{-1}(x_i)$ the stabilizer is $G_x=\Z/\fe_i\Z$, where $\fe_i$ is the ramification index. 
We identify $G_x$ with $\Z/\fe_i\Z$ via the generator $\sigma_i=\sigma^{m/\fe_i}$.
We write $\gamma_i=\sigma_i^{k_i}$, where $k_i$ and $\fe_i$ are co-prime.
The action of $\Z_m$ induces a character $\chi_x$ of $G_x$ via its action on the tangent space at $x$, which is related to the holonomy data via
$$\chi_x(\sigma_i)=\zeta^{\nu_i \frac{m}{\fe_i}},$$
where $\nu_i$ is the inverse of $k_i$ in $\Z_{\fe_i}$.
The algebra $p_*\cO_C$ is a locally free sheaf of rank $m$ with a $\Z_m$-action, and admits a decomposition into isotypical factor $$p_*\cO_C = \bigoplus_{j=0}^{m-1} T_j,$$ where $T_j$ is the subsheaf of sections $s$ such that $\sigma^*s = \zeta^js$.

\begin{lem}\label{lem:multracinetrivial}
There is a canonical morphism $T_1^{\otimes m}\to \cO_D$, which is an isomorphism.
The multiplicity of $T_1$ at $x_i$ (i.e. the representation $\bmu_{\fe_i}\to \G_m$ given by $L_{|x_i}$) is $-\nu_i/\fe_{i}$.
\end{lem}
Thus, we obtain an isomorphism $$\overline{\cM}_{g,n}(B\Z_m,\ug)\simeq \overline{\cM}^{m,0}_{g,\underline{b}},$$
where $\ub$ is given by $b_i = -\nu_i \frac{m}{\fe_i}$.
The description of $\overline{\cM}_{g,n}(B\Z_m)$ in terms of roots of the trivial bundle is more convenient in the next section.

\subsection{The spine CohFT}

\label{subsection:spinecohft}

Let $\ua \in (\Z/rm\Z)^n$ and $\ub\in (\Z/m\Z)^n$ be multi-indices.
Define the space $\overline{\cM}^{rm}_{0}(B\Z_m,\ua,\ub)$ to be the stack with objects $\overline{\cM}^{rm}_{0}(B\Z_m,\ua,\ub)(S) =\left\{(\cC,\overline{\cC},T,\alpha_{1},\alpha_{2}) \right\} $ where
\begin{itemize}
	\item $\cC$ is a stable twisted curve over $S$,
	\item $\overline{\cL}$ and $T$ are line bundle over $\cC$ having multiplicity $\frac{\ua}{rm}$ and $\frac{\ub}{m}$, such that $\cL \oplus T$ is representable,
	\item $\alpha_{1}:\overline{\cL}^{\otimes rm}\to \omega_{\log}$ and $\alpha_{2}: T^{\otimes m} \to \cO_{\cC}$ are isomorphisms.
\end{itemize}

Thus, the universal curve carries two universal line bundles $\overline{\cL}_{\ua}$, and $T_{\ub}$.
For a couple of multi-indices $\ua,\ub$, we define the subsets 
\begin{align*}
&	S^{0}= \left\{ i | a_{i}=0 \wedge b_{i}=0\right\}
&	S^{\times}_{j}= \left\{ i | \mult_{x_{i}}(\overline{\cL}\otimes T^{j}=0 \wedge b_{i}\in \Z_{m}^{\times} \right\},
\end{align*}
and the associated divisors in the universal curve
\begin{align*}
&	E^{0}= \bigsqcup_{i\in S^{0}} x_{i},&
	E^{\times}_{j}= \bigsqcup_{i\in S^{\times}_{j}} x_{i}.
\end{align*}

Let $\left\{\cA_{\xi,j}| j\in \Z_{m}, \xi \in \bmu_{rm}\right\}$ be invertible multiplicative classes with $\cA_{\xi,j}=\exp( \sum_{d\geq 0} w^{\xi,j}_d \ch_d)$.
We introduce the following classes on $\overline{\cM}^{rm}_{0}(B\Z_{m},\ua,\ub)$
\begin{equation}
\Lambda^{\spine}_{\xi}(\ua,\ub)_{0,n} = r m^2\prod_{j=0}^{m-1} \cA_{\xi,j}\left( R\pi_*\left( \overline{\cL}_{\ua}\otimes T^{j}_{\ub}(-E^{0}-E^{\times}_{j})\right) \right).
\end{equation}

\begin{defn}[Spine state space]
	The state space of the spine CohFT is 
	$$V_{\spine} = \C[\bmu_{rm}]\otimes \C\left[\widehat{\bmu}_{rm}\right]\otimes \C\left[\Z_m\right]\otimes \C\llbracket s \rrbracket.$$
	We fix the basis $\{\phi_{a,\xi} \otimes [b]\}_{a,\xi,b}$ of $V_{\spine}$ as a free $\C\llbracket s \rrbracket$-module.
	The pairing is given by
\begin{equation}
		\scal{\phi_{a,\xi} \otimes [b];\phi_{a',\xi'} \otimes [b']} = \left\{ \begin{matrix} 
			\delta_{a,-a'}\delta_{\xi,\xi'}\delta_{b,-b'} &\textrm{ if } a\neq 0[r]\\
			\delta_{a,-a'}\delta_{\xi,\xi'}\delta_{b,-b'}\prod_{j\in\Z_{m}}\cA_{\xi,j}^{-1}(\cO) &\textrm{ if } a=0 \textrm{ and } b=0,\\
		\delta_{a,-a'}\delta_{\xi,\xi'}\delta_{b,-b'}\cA_{\xi,j}^{-1}(\cO) &\textrm{ if } \frac{a}{rm}+\frac{jb}{m}=0[1] \textrm{ and } b\in\Z_{m}^{\times}.
		\end{matrix} \right.
\end{equation}
\end{defn}

\begin{rem}
As in the previous section, $V_{\mathrm{spine}}$ is the direct sum of $rm$ different copies of $\C[\Z_{rm}] \otimes \C[\Z_{m}]\otimes \C\llbracket s \rrbracket$ indexed by $\hat{\bmu}_{rm}$.
\end{rem}

\begin{prop}
	The projection to $\overline{\cM}_{0,n}$ of the classes $\Lambda^{\spine}_{\xi}$ form the genus-$0$ part of a CohFT over $V_{\spine}$.
\end{prop}

\begin{defn}
	The CohFT defined above is called the spine CohFT, and its associated Lagrangian cone is denoted by $L^{\spine}$.
\end{defn}

We now apply Chiodo--Zvonkine's theorem to compute the cone $L^{\spine}$ in the relevant sectors.
\begin{prop}
There exists a linear operator $\Box$ such that the cone $L^{\spine}$ is related to the untwisted cone via
$$L^\spine = \Box L^\un.$$
Moreover, we have

\begin{equation}
	\Box \cdot \phi_{a ,\xi}\otimes [b] =  \prod_{j\in \Z_m} \exp\left( \sum_{d\geq 0} w_d^{\xi,j} \frac{z^d}{(d+1)!} \widetilde{B}_{d+1}\left( \frac{\fa+rl}{rm} + \frac{jb}{m} \right)\right)\phi_{a,\xi}\otimes [b],
\end{equation}
where $\widetilde{B_{d}}$ is the restriction of the Bernoulli polynomial to $]0;1]$ (taken $\Z$-periodically).

\end{prop}

\begin{proof}
This is a straightforward extension of \cite{chiodoTwistedGromovWittenRspin2007}, theorem 1.2.2.
The restriction of the Bernoulli polynomial to $]0;1]$ instead of $[0;1[$ comes from the $-E$ twisting at the broad points in the definition of the invariants.
\end{proof}

We now apply the previous result to the classes encountered in Lefschetz formula (see \Cref{adelic}).
Let $\cA_{\xi,j}$ be the classes with values in $\C\llbracket s \rrbracket$ defined by
\begin{equation}
\cA_{\xi,j}(E)=\ch\left( \lambda_{-s\xi e^{\frac{-2i\pi j}{m}}}E\right)^{-(N+1)}.
\end{equation}
They correspond to the power series
$$w_d^{\xi,j} =(N+1) \sum_{k\geq 1} \frac{\left( se^{\frac{-2i\pi j}{m}}\right)^k\xi^{k} k^d}{k}.$$

Let $\xi_0=\exp(2i\pi k_{0}/rm)\in \bmu_{rm}$ be a root of unity such that $\gcd(k_{0},m)=1$, and let $\nu_{0}$ be the inverse of $k_{0}$ modulo $m$.
The restrictions of $\Box$ to the sectors $b=0$, and $b=\nu_{0}$ are denoted by $\Box_0$ and $\Box_{\xi_0}$ respectively, and are given by
\begin{align*}
	\MoveEqLeft[3]	\Box_0 \cdot \phi_{a+1,\xi}\\{}  &= \exp\left( (N+1)\sum_{k,d,j} \frac{\left( se^{\frac{-2i\pi j}{m}}\right)^k\xi^{k} k^d}{k(d+1)!}z^d B_{d+1}\left(\frac{a+1}{rm}\right)\right)\phi_{a+1,\xi}\\
				     &=  \exp\left((N+1) \sum_{k,d}m \frac{ (s\xi)^{mk} (km)^d}{mk(d+1)!}z^d B_{d+1}\left(\frac{a+1}{rm}\right)\right)\phi_{a+1,\xi}\\
&= \exp\left((N+1) \sum_{k} \frac{s^{mk}\xi^{mk}}{k}  \frac{e^{\frac{(a+1)kmz}{rm}}}{e^{kmz}-1}\right)\phi_{a+1,\xi}\\
&=\exp\left((N+1) \sum_{k} \frac{s^{mk}\xi^{mk}}{k}  \frac{q^{\frac{(a+1)k}{r}}}{q^{km}-1}\right)\phi_{a+1,\xi}  \textrm{ for } 0\leq a <rm.
\end{align*}
Let $\widetilde{B}_d(x)$ denote the Bernoulli polynomial restricted to $]0;1]$, and expanded $\Z$-periodically. 
For $a\notin m\Z_{rm}$, let us write $a = \fa + rl$, with $1\leq \fa \leq r-1$. Then we have
\begin{align*}
	\MoveEqLeft[3] \Box_{\xi_0} \cdot \phi_{a,\xi} \\ {} &=   \exp\left((N+1) \sum_{k,d,j} \frac{s^k \xi^{k} e^{\frac{-2i\pi kj}{m}}(kz)^d}{k(d+1)!} \widetilde{B}_{d+1}\left( \frac{a}{rm} + \frac{\nu_0j}{m}\right) \right) \phi_{a,\xi} \\
&=\exp\left((N+1) \sum_{k,d,j} \frac{ s^k\xi^{k}\xi_0^{-rjk} (kz)^d}{k(d+1)!} \widetilde{B}_{d+1}\left(\frac{a}{rm}+\frac{j}{m}\right)\right)\phi_{a,\xi}\\
&=\exp\left((N+1) \sum_{k,d,j} \frac{ s^k\xi^{k}\xi_0^{-rjk}\xi_0^{rlk} (kz)^d}{k(d+1)!} {B}_{d+1}\left(\frac{\fa }{rm}+\frac{j}{m}\right)\right)\phi_{a,\xi}\\
&= \exp\left( (N+1) \sum_{k,j} \frac{ s^k\xi^{k}\xi_0^{-rjk}\xi_0^{rlk}}{k} \frac{q^{k\frac{\fa+rj}{rm}}}{q^{k}-1} \right) \phi_{a,\xi}\\
&= \exp\left((N+1) \sum_{k} \frac{ s^k\xi^{k}\xi_0^{rlk}}{k} \frac{q^{k\frac{\fa}{rm}}}{q^{k/m}\xi_0^{-rk}-1} \right) \phi_{a,\xi}.
\end{align*}

Similarly, for $a=rl$ we have
\begin{align*}
	\MoveEqLeft[3] \Box_{\xi_0} \cdot \phi_{rl,\xi} \\{} &=   \exp\left( (N+1)\sum_{j=0}^{m-1}\sum_{k,d} \frac{s^k \xi^{k} \xi_0^{-rkj}}{k(d+1)!}z^d \widetilde{B}_{d+1}\left(\frac{l+j}{m}\right) \right) \phi_{rl,\xi}\\
							     &=\exp\left((N+1)\sum_{j=1}^m\sum_{k,d} \frac{ s^k\xi^{k}\xi_0^{rk(l-j)} (kz)^d}{k(d+1)!} B_{d+1}\left(\frac{j}{m}\right)\right)\phi_{rl,\xi}\\
							     &= \exp\left((N+1) \sum_{j=1}^m\sum_{k \geq 1} \frac{ s^k\xi^{k}\xi_0^{rk(l-j)}}{k}\frac{q^{kj/m}}{q^{-k}-1}\right)\phi_{rl,\xi}\\
							     &= \exp\left((N+1) \sum_{k \geq 1} \frac{ s^k\xi^{k}\xi_0^{rkl}}{k} \frac{q^{k/m}}{\xi_0^{-rk}q^{k/m}-1}\right)\phi_{rl,\xi}.
\end{align*}

\begin{rem}
	We chose to express the operator $\Box$ in terms of an $rm$th root of unity $\xi_{0}$ to make the connection with quantum K-theory easier in the next section, where $\xi_{0}$ will be the trace of an automorphism on the first cotangent line bundle.
	Here, the resulting operator only depends on the $m$th root $\xi_{0}^{r}$, or $\nu_{0}\in \Z_{m}$.
\end{rem}


\subsection{Spine invariants}
Similarly to \Cref{subsection:fakeinvariants}, we multiply the classes $\Lambda^\spine_\xi(\ua,\ub)$ by classes $\cB$ and $\cC$ to define the spine invariants.
For our application, we are only interested in two kinds of sectors: $b=0$ and $b\in \Z_{m}^{\times}$, so we only introduce the twists there.  
It is a consequence of \cite{tonitaTwistedOrbifoldGromovWitten2014} that the $J$-function in these sectors does not depend on the twists elsewhere.

The singular locus $Z\subset \overline{\cC}$ admits a decomposition
$$Z = \bigsqcup_{\substack{\overline{\fe} |r \\ \ff|m}} Z_{\overline{\fe},\ff},$$
where $Z_{\overline{\fe},\ff}$ is made of the nodes such that $\mult(\overline{\cL})$ has order $\overline{\fe}$, and $\mult(T)$ has order $\ff$.
We also denote by $\fe(a)$ the order of $a\mod r$.

\begin{defn}[Spine invariants]
\label{defn:spine_class}
We define the classes 
\begin{align*}
\cB_{n} &= \td^{-1}(\pi_*\omega_{\log}^{-1}-1) \prod_{k=1}^{m-1} \td_{\zeta^k}\left( \pi_*\left( 1-\omega_{\log}^{-1}\otimes T^k\right) \right),\\
\cC_{n}^{\bar{\fe},1} &= \td^{\vee}(-\pi_*\cO_{Z_{\bar{\fe},1}})\prod_{k=1}^{m-1}\td_{\zeta^k}^{\vee}(-\pi_{*}(\cO_{Z_{\bar{\fe},1}})),\\
\cC_{n}^{\bar{\fe},m} &= \td^{-\otimes{\lcm(\bar{\fe},m)/m\fe}}(-\pi_{*}\cO_{Z_{\bar{\fe},m}})
\end{align*}
where $\zeta = \exp( \frac{2i \pi }{m})$, and 
\begin{align*}
\td_\zeta(L) &= \frac{1}{1-\zeta \exp(-c_1(L))},\\
\td^{\otimes l}(L)&= \frac{lc_{1}(L)}{1-e^{-lc_{1}(L}}.
\end{align*}

The spine invariants are defined by
\begin{equation}
	\scal{\phi_{a_1,\xi}\otimes [b_1]\cL_1^{k_1},\hdots,\phi_{a_n,\xi}\otimes [b_n]\cL_n^{k_n}}^{\spine}_{0,n} = \int_{\overline{\cM}^{rm}(B\Z_m,\ua,\ub)}\Lambda_\xi^{\spine}(\ua,\ub)\cB_n\cC_n^{\bar{\fe},0}\cC_{n}^{\bar{\fe},m}.
\end{equation}
\end{defn}

\begin{rem} Let us sketch the origin of the twisting classes $\cB$ and $\cC$ (see also \cite[Section 8]{giventalHirzebruchRiemannRoch2011}).
	The moduli stack of spines (see \Cref{defn:spine}) parametrizes $r$-spin curves $\cC$ together with a $\Z_{m}$-action permuting the marked points.
	The twisting classes above arise in Lefschetz formula as the class
	\begin{equation}
\label{equ:Lefschetzclass}
		\frac{\td(T_{\cM^{\spine}})}{\ch \circ \Tr(\lambda_{-1}\cN^{\vee})},
	\end{equation} where $\Tr$ denotes the trace bundle in K-theory, and $\cN$ is the normal bundle to the morphism $\cM^{\spine}\to \overline{\cM}^{r}_{0,n}$.
	The pullback of the tangent sheaf of $\overline{\cM}^{r}_{0,n}$ is equipped with a $\Z_{m}$-action, and admits a decomposition into a fixed and moving part.
	The subsheaf of invariant sections is $\cT_{\cM^{\spine}}$, and the moving part is the normal bundle $\cN$.
	Let us describe this decomposition.
	The moduli stack of spines in equipped with $2$ universal curves $\cC$ and $\cD$, which fit in the following diagram
	\begin{equation}
		\begin{tikzcd}
			\cC \arrow[r,"p"] \arrow[d,"\pi"] & \cD=\cC/\Z_{m} \arrow[dl,"\overline{\pi}"]\\
			\cM^{\spine}.
		\end{tikzcd}
	\end{equation}
	The tangent sheaf of $\overline{\cM}^{r}$ is 
\begin{align*}
	\cT_{\overline{\cM}^{r}_{0,n}} &= -R\pi_{*}(\omega_{\log}^{\vee}) -(\pi_{*\cO_{Z}})^{\vee}\\
				       &= -R\overline{\pi}_{*}(p_{*}\omega_{\log}^{\vee}) - \overline{\pi}_{*}(p_{*}\cO_{Z})\\
				       &= -\bigoplus_{k=0}^{m-1} R\overline{\pi}_{*}(\omega_{\log}^{\vee}\otimes T^{k}) - \overline{\pi}_{*}(p_{*}\cO_{Z}).
\end{align*}
In the last equation, the first term is decomposed into eigensheaves; so its contribution to the class \eqref{equ:Lefschetzclass} is $$\td\left(R\overline{\pi}_{*}\omega_{\log}^{\vee}\right)^{-1}\prod_{k=1}^{m-1}\td_{\zeta^{k}}\left(-\overline{\pi}_{*}\omega_{\log}^{\vee}\otimes T^{\otimes -k}\right).$$
For the second term, we may assume that the ramification index of $p$ at every node is either $1$ or $m$.
Then we obtain 
$$p_{*}\cO_{Z} = \bigoplus_{\fe}\cO_{Z_{\fe,1}}[\Z/m\Z] \oplus \cO_{Z_{\fe,m}}.$$
The action on the first term is the regular representation of $\Z_{m}$, and the action on the second term is trivial.
Thus, its contribution to \ref{equ:Lefschetzclass} is 
$$\td^{\vee}(-\overline{\pi}_{*}\cO_{Z_{\fe,1}})\td^{\vee}(-\overline{\pi}_{*}\cO_{Z_{\fe,m}}) \prod_{k=1}^{m-1}\td^{\vee}_{\zeta^{k}}(-\pi_{*}\cO_{Z_{\fe,1}}).$$
\end{rem}

\begin{rem}
	The exponent in the definition of $\cC^{\bar{\fe},m}_{0,n}$ comes from a slight difference between $\cM^{\spine}$ and $\overline{\cM}^{rm,m}(B\Z_{m})$: the nodes of the quotient curve have stabilizers of order $m\fe$ in the first case, and $\lcm(\bar{\fe}(a),m)$ in the second case (see \Cref{defn:spine}).
	Thus, a correction is needed to match the integrals of the spine contribution (see \Cref{subsection:Spine contribution}).
\end{rem}
\begin{prop}[{\cite[cor. 6.2, 6.3]{tonitaTwistedOrbifoldGromovWitten2014}, \cite[Section 7]{giventalHirzebruchRiemannRoch2011}}]
	\label{prop:polarisation spine}
The $\cB$ twist changes the dilaton shift to $1-q^m$.
In the $0$ sector, the $\cC$ twist changes to polarization to 
\begin{equation}
	\cK^{\spine,0}_{-} = \mathrm{Span}\left\{ \frac{q^{mk/\bar{\fe}(a)}}{(1-q^{m/\bar{\fe}(a)})^{k+1} }\phi_{a,\xi} \otimes [0] | k\in \N\right\}.
\end{equation}
For $b \in \Z_m^{\times}$, the polarization in the $b$ sector is given by 
\begin{equation}
	\cK^{\spine,b}_{-} = \mathrm{Span}\left\{ \frac{q^{k/m\fe(a)}}{(1-q^{1/m\fe(a)})^{k+1}}\phi_{a,\xi} \otimes [b] | k\in \N\right\},
\end{equation}
\end{prop}

Thus, the $J$-function of the spine invariants is defined by 
\begin{align*}
	\MoveEqLeft[3]	J^{\spine}(t)= 1-q^m + t \\
	{} &+ \sum_{n,a,\xi} \frac{m \phi^{a,\xi}\otimes [0]}{\bar{\fe}(a)n!} \corr{\frac{\phi_{a,\xi} \otimes [0]}{1-q^{m/\bar{\fe}}\cL_{0}^{m}},t(\cL_{1}),\hdots,t(\cL_{n})}^{\spine}_{0,n+1} \\
	   &+\sum_{\substack{a\in \Z_{rm} \\ b \in \Z_{m}^{\times} \\n \geq  2}} \frac{\phi^{a,\xi} \otimes [b]}{m\fe(a) n!} \corr{\frac{\phi_{a,\xi} \otimes [-b]}{1-q^{\frac{1}{m\fe}}\cL_0^{\frac{\lcm(\bar{\fe},m)}{m\fe}}},t(\cL_1),\hdots , t(\cL_n)}^\spine_{0,n+1}\\
	   &+ \textrm{ other sectors.}
\end{align*}
Recall that $\fe(a)$ is the order of $a$ in $\Z_{r}$, while $\bar{\fe}(a)$ is the order of $a$ in $\Z_{rm}$.
The spine invariants satisfy a natural symmetry, arising from the cyclic permutation of the line bundles $\overline{\cL}\otimes T_j$ and from our choice of $\cA_{\xi,j}$.

\begin{lem}
\label{lem:invariancespine}
We have
\begin{equation}
\Lambda^\spine_{\xi}(\ua,\ub) = \Lambda^\spine_{\xi e^{\frac{-2i\pi}{m}}}(\ua + r\ub,\ub).
\end{equation}
In particular, if $T$ is an element of $\cK^\fake_+$, $\zeta\in \bmu_{rm}$, and  $a'\in \Z_{rm}$, then the following correlator does not depend on $k\in \Z_m$
\begin{equation}
\scal{\phi_{a+rk,\xi\xi_0^{a+rk}}\otimes [-\nu_0], \sum_{k'=0}^{m-1}\phi_{a'+rk',\zeta \xi_0^{-a'-rk'}}\otimes [\nu_0], \Psi^m(\Phi_0(T)),\hdots,\Psi^m(\Phi_0(T))}_{0,n+2}^{\spine},
\end{equation}
where the Adams operation acts by $\Psi^{m}(\phi_{a,\xi}\otimes [0]) = \sum_{\zeta^{m}=\xi} \phi_{a,\zeta}\otimes [0]$.
\end{lem}

\subsection{Comparison with the fake theory}
Recall that the Adams operation in $K$-theory are the ring morphisms $\Psi^{m}:K^{0}(X)\to K^{0}(X)$ such that for any line bundle $L$, we have $$\Psi^{m}(L)=L^{\otimes m}.$$
The Adams operations extend to the state space $\cK^{\spine}$ by the formula
$$\Psi^{m}\left(s^{k} q^{j/rm}\phi_{a,\xi}\otimes [b]\right) = s^{km}q^{j/r} \sum_{\zeta^{m}=\xi}\phi_{a,\zeta}\otimes [b].$$

We now show that the $0$-sector of the spine cone $L^{\spine}$ contains $\Psi^m(L^\fake)$.

\begin{lem}
\label{lem:adamsuntwited}
The untwisted cone $L^\un \subset \cK^\spine$ is stable by the transformation $q\mapsto q^m$.
\end{lem}
\begin{proof}
The cohomological cone of a point is $$L_{\pt} = z\bigcup_{\tau\in \C} \exp\left(\frac{\tau}{z}\right)\cH_+,$$
which is obviously invariant by the transformation $z\mapsto mz$.
The result follows for the untwisted cone.
\end{proof}

\begin{lem}
\label{lem:morphismesymplectique}
Let $\Phi_0: \cK^\fake \to \cK^\spine$ be the morphism $\phi_{a}\otimes e_\xi \mapsto \phi_{ma,\xi} \otimes [0]$.
Then, $\Psi^m \circ \Phi_0$ is a morphism of polarized symplectic spaces.
\end{lem}

\begin{proof}
This follows from \Cref{prop:polarisation spine}.
\end{proof}

\begin{prop}
\label{prop:psim(delta)=box0}
In $\cK^\spine$ we have that 
\begin{equation}
\Psi^m\left( \Delta \Phi_0L^\un\right)  = \Psi^m\left( \Phi_0\Delta L^\un \right)\subset \Box_0 L^\un.
\end{equation}
In particular for any element $t$ of $\cK^\fake_+$, we have
$$J^\spine(\Psi^m \Phi_0 t ) = \Psi^m \Phi_0 J^{\fake}(t).$$
\end{prop}

\begin{proof}
The first assertion comes from a direct computation, together with \Cref{lem:adamsuntwited}.
Indeed, if $f\in \cK^{\fake}$ is an element of $L^{\un}$, then $\Psi^{m}\Phi_{0}f \in L^{\un}$.
Moreover, we compute that for any $\xi \in \bmu_{r}$ and $a\in \{0,\hdots, r-1\}$
\begin{align*}
	\Psi^{m}\left( \Delta \Phi_{0}(\phi_{a+1,\xi})\right) &= \sum_{\zeta^{m}=\xi} \exp\left( r\sum_{k\geq 1} \frac{s^{km}\zeta^{mk}}{k} \frac{q^{mk(a+1)/r}}{q^{mk}-1}\right) \phi_{m(a+1),\zeta}\\
								       &= \Box_{0} \left( \sum_{\zeta^{m}=\xi} \phi_{m(a+1),\zeta}\right)\\
								       &= \Box_{0}\Psi^{m}\Phi_{0}(\phi_{a+1,\xi}).
\end{align*}
So in particular we have that $\Psi^{m}\Phi_{0}\Delta f = \Box_{0} \Psi^{m}\Phi_{0}f \in \Box_{0}L^{\un}$.

The second assertion is deduced form the first one, together with the change in dilaton shift from $1-q$ to $1-q^m$, and \Cref{lem:morphismesymplectique}.
\end{proof}

\begin{cor}
\label{prop:twistedfakepsim}
Let $J_{(1)}(t)$ be the image in $\cK^\fake$ of $J(t)$, and let $T=\left[J_{(1)}\right]_+-1+q$, where $[\cdot]_+$ denotes the projection on $\cK^\fake_+$.
Then, $$J^\spine\left( \Psi^m\Phi_0 T\right) = \Psi^m \Phi_0\left( J^\fake(T)\right) = \Psi^m \Phi_0 J_{(1)}(t).$$
\end{cor}

\section{Adelic characterization}
\label{adelic}
In this section, we use the Grothendieck--Riemann--Roch theorem \cite{toenTheoremesRiemannRoch1999} to compute the $J$-function in terms of the fake and spine theories (see \Cref{section:spinecohft} for the definition of the spine CohFT).
We first recall that for  equivariant sheaves, this theorem takes the form of a Lefschetz fixed-point formula,  and we deduce a natural $\bmu_r$-action on the loop space making the $J$-function equivariant.
This action can be interpreted as an automorphism of the sections of the universal curve.
Then, we use Lefschetz formula to compute the expansion of the $J$-function at each root of unity.
The results are expressed in the adelic characterization theorem, which characterizes values of the $J$-function in terms of their expansions at each root of unity.

\subsection{Lefschetz formula and the $\bmu_r$-action}
\label{subsection:mu_raction}

Let $X$ be a proper smooth Deligne--Mumford stack over $\C$, let $h$ be an automorphism of $X$ of finite order, and let $F$ be an equivariant coherent sheaf.
Then Lefschetz formula \cite{toenTheoremesRiemannRoch1999} reads
\begin{equation}
\tr_h\left( H^*\left(X,f\right) \right) = \int_{X^h} \ch\left( \frac{ \Tr_h(F)}{\Tr_g \lambda_{-1}\cN^\vee} \right) \Td(T_{X^h}),
\end{equation}
where $X^h$ is the fixed-point stack, and $\cN$ is the normal bundle to the morphism $X^h\to X$.
Finally, $\Tr_h (F)$ is obtained by decomposing $F_{|X^h}$ into isotypical factors $F=\bigoplus_\lambda F^\lambda$, and multiplying each factor by $\lambda$ 
\begin{equation}
\Tr_h(F)= \sum_\lambda \lambda F^\lambda \in K^0(X)\otimes \C.
\end{equation}

This formula suggests that the potential $\cF$, and the $J$-function are equivariant with respect to a natural action of $\bmu_r$ on $\cK$.
Indeed, for all $h\in S_n$, objects of the fixed-point stack $(\widetilde{\cM}^r_{0,n+1})^h$ are given by the data $(C,\cL,(\sigma_i)_{i=0}^n,g,\phi,(\eta_i)_{i=0}^n)$, where
\begin{itemize}
\item $(C,\cL,(\sigma_i)_{i=0}^n)$ is an object of $\widetilde{\cM}^r_{0,n+1}$ over $S$,
\item $g$ is an automorphism of $C$,
\item $\phi : g^*\cL \to \cL$ is an isomorphism compatible with the spin structure,
\item $\eta_i:g\circ \sigma_i\to \sigma_{h(i)}$ is a $2$-isomorphism.
\end{itemize}

We let $\bmu_r$ act on $( \widetilde{\cM}^r_{0,n+1})^h$ by changing the $2$-isomorphisms $\eta_0,\hdots, \eta_{n}$.
Explicitly, if $\zeta$ is an element of $\bmu_r$, an object $((C,\cL,(\sigma_i)_{i=0}^n,g,\phi,(\eta_i)_{i=0}^n)$ is sent to 
\begin{equation}
\zeta \cdot (C,\cL,(\sigma_i)_{i=0}^n,g,\phi,(\eta_i)_{i=0}^n) = (C,\cL,(\sigma_i)_{i=0}^n,g,\phi,\zeta^{\frac{r}{\fe_i}}\circ \eta_i),
\end{equation}
where $\fe_i$ denotes the cardinality of the stabilizer at the $i$th marked point. 

This action has a natural analogue on the state space.
\begin{defn}
We define a $\bmu_r$-action on $\cK$ by 

\begin{align}
	\zeta\cdot\left( f( q^{1/r}) \phi_{a,l}\right)& = f(\zeta^{-1} q^{1/r})\zeta^{al}\phi_{a,l}.
\end{align}

\end{defn}
\begin{rem}
	In the idempotent basis, the action becomes
	\begin{equation}
		\zeta \cdot (f\phi_{a}e_{\xi}) = f(\zeta^{-1} q^{1/r}) \phi_{a} e_{\xi\zeta^{a}}.
	\end{equation}
\end{rem}
\begin{rem}
	\label{rem:action}
	The action in the previous definition is designed to be compatible with the $\bmu_{r}$-action on the fixed-point stack $\left( \widetilde{\cM}_{0,n+1}^{r}\right)^{h}$ in the following sense:
	\begin{equation}
		\zeta^{*}\left(\ch\circ \Tr \left( \frac{\Lambda_{n}(s)\bigotimes_{i=1}^{n} \ev_{i}^{*}(t(\cL_{i}))}{\lambda_{-1}\cN^{\vee}} \right) \right) =\ch \circ \Tr \left( \frac{\Lambda_{n}(s) \bigotimes_{i=1}^{n} \ev_{i}^{*}(\zeta\cdot t (\cL_{i})}{\lambda_{-1}\cN^{\vee}}\right). 
	\end{equation}
\end{rem}

\begin{prop}
The genus-0 potential is $\bmu_r$-invariant, and the $J$-function is $\bmu_r$-equivariant:
\begin{equation}
\zeta \cdot J(t)= J(\zeta \cdot t).
\end{equation} 
\end{prop}

\begin{proof}
The statement about the potential is a consequence of the previous remark, together with the Lefschetz formula.
To prove the second statement, we write 
\begin{align*}
	\zeta \cdot \left( J(t) \right)&= 1-q + \zeta \cdot t + \sum_{\substack{a\in \Z_{r} \\ \xi \in \bmu_{r}}} \phi^{a,\xi\zeta^{-a}}\scal{\frac{\phi_{a,\xi}}{1-q^{1/\fe(a)}\zeta^{-r/\fe(a)}\cL_{0}},t,\hdots,t }^{S_{n}}_{0,n+1}\\
				       &= 1-q + \zeta\cdot t+\sum \phi^{a,\xi\zeta^{-a}}\scal{\frac{\phi_{a,\xi}}{1-q^{1/\fe(a)}\zeta^{-r/\fe(a)}\cL_{0}},\zeta \cdot t,\hdots,\zeta \cdot t }^{S_{n}}_{0,n+1}\\
				       &= 1-q +\zeta\cdot t+\sum \phi^{a,\xi\zeta^{-a}}\scal{\zeta \cdot \frac{\phi_{a,\xi\zeta^{-a}}}{1-q^{1/\fe(a)}\cL_{0}},\zeta \cdot t,\hdots,\zeta \cdot t }^{S_{n}}_{0,n+1}\\
				       &=1-q + \zeta\cdot t+\sum \phi^{a,\xi\zeta^{-a}}\scal{\frac{\phi_{a,\xi\zeta^{-a}}}{1-q^{1/\fe(a)}\cL_{0}},\zeta \cdot t,\hdots,\zeta \cdot t }^{S_{n}}_{0,n+1}\\
				       &= J(\zeta \cdot t).
\end{align*}
\end{proof}

Let $t$ be an element of $\cK_{+}$, and let $\underline{t} =\frac{1}{r}\sum_{\zeta\in\bmu_{r}} \zeta \cdot t$ denote its projection onto the subspace of invariant vectors. 
As a consequence of the previous proposition, we have that $$\cF_{0}(t)= \cF_{0}(\underline{t}).$$
This also follows from \Cref{rem:vanishing}, because $\bmu_{r}$-invariant elements of $\cK_{+}$ are precisely those elements which satisfy the non-vanishing condition.

\begin{cor}
For all $t\in \cK_+$, the projection of $\left[J(t)\right]_-$ of $J(t)$ to $\cK_-$ parallel to $\cK_+$ is a $\bmu_{r}$-invariant point.
\end{cor}
\begin{proof}
	By the preceeding proposition we have that $[J(t)]_{-} = [J(\underline{t})]_{-}$.
	But $\underline{t}$ is $\bmu_{r}$-invariant, so $J(\underline{t})$ is a $\bmu_{r}$-invariant point.
	\end{proof}
\subsection{Adelic characterization}
Following \cite{giventalHirzebruchRiemannRoch2011} and \cite{giventalPermutationEquivariantQuantumKtheorya}, we apply Lefschetz, formula to the $J$-function to find recursion relations.
More precisely, we compute the expansion of the $J$-function at each root of unity $\xi_0$, and show that it corresponds to the fake and spine theories.

\begin{defn}
\label{defn:morphismePhi}
Let $\xi_0$ be a root of unity, and let $m$ be the order of $\xi_{0}^r$.
Thus we can write $\xi_{0} = e^{\frac{2i \pi k_{0}}{rm}}$, where $k_{0}$ is prime to $m$.
Finally, let $\cK_{rm}$ be the symplectic space defined in \Cref{defn:krm}.

We define linear maps $\Phi_{0},\Phi_{\xi_{0}} : \cK \to \cK_{rm}$ by
\begin{align*}
	\Phi_{0} : \cK &\to \cK_{rm} \\
	f(q^{1/r})\phi_{a,\xi} &\mapsto f(q^{1/r})\phi_{ma,\xi}, 
\end{align*}
and 
\begin{align*}
	\Phi_{\xi_{0}} : \cK &\to \cK_{rm} \\
	f(q^{1/r}) \phi_{a,\xi} &\mapsto f(\xi_{0}^{-1}q^{1/r}) \sum_{l=0}^{m-1}\phi_{a+rl,\xi \xi_{0}^{-a-rl}}.
\end{align*}

We think of $\Phi_0$ (resp. $\Phi_{\xi_0}$) as an embedding of $\cK$ in the $0$ sector (resp. the $k_0^{-1}$ sector) of the spine CohFT. 
\end{defn}
\begin{thm}[Adelic characterization]
\label{thm:adelic_char}
Let $f(q^{1/r})$ be a $\bmu_r$-invariant element of $\cK$ such that $f\in 1-q + \scrI\cK$.
Then $f$ lies is the image of the $J$-function $L^{K}_{\mathrm{FJRW}}$ if an only if 
\begin{itemize}
\item The poles of $f$ belong to $\bmu_\infty \cup \{0,\infty\}$,
\item the expansion $f_{(1)}$ at $q^{1/r}=1$ belongs to the fake cone $L^\fake$,
\item For all $\xi_0\in \bmu_\infty$ such that $\xi_0^r$ has order $m$,
we have $$\Phi_{\xi_0} \left(f\right) (q^{\frac{1}{rm}})\in \Box_{\xi_0}\Delta^{-1}\cT L^\fake, $$ where $\cT L^\fake$ is the tangent space at the point $\Phi_0f \in L^\fake$, $\Delta$ is the operator of the fake theory (\Cref{prop:Delta}), and for $a\in \{1,\hdots,r-1\}$, the operator $\Box_{\xi_0}$ is defined by
\begin{align*}
\Box_{\xi_0}(\phi_{a+rl,\xi} ) &= \exp\left( (N+1)\sum_{k\geq 1}\frac{s^k\xi^{k}\xi_0^{rlk}}{k} \frac{q^{k\frac{a}{rm}}}{\xi_0^{-rk}q^{k/m}-1} \right) \phi_{a+rl,\xi}\\
\Box_{\xi_0}(\phi_{rl,\xi} ) &= \exp\left( (N+1)\sum_{k\geq 1}\frac{s^k\xi^{k}\xi_0^{rlk}}{k} \frac{q^{\frac{k}{m}}}{\xi_0^{-rk}q^{k/m}-1} \right) \phi_{rl,e_\xi}.
\end{align*}
\end{itemize}
\end{thm}

\begin{rem}
We showed in \Cref{subsection:mu_raction} that the $J$-function is $\bmu_r$-equivariant, and that for all $t\in \cK_+$, we have $$\left[J(t)\right]_- =\left[J(\underline{t})\right]_- ,$$
where $\underline{t} = \frac{1}{r} \sum_{\xi \in \bmu_r} \xi \cdot t$ is the projection of $t$ on the subspace of invariant elements.
Thus, the theorem above characterizes all possible values of the $J$-function.
\end{rem}

The rest of this section is devoted to the proof of \Cref{thm:adelic_char}.
We first show that these 3 conditions are necessary. 
The first item is obviously necessary, and the other two follow from Lefschetz formula, which allows us to compute the expansion of $J(t)$ at $q^{1/r} =1$, and $q^{1/r} = \xi_0$.

\begin{defn}
Let $\xi_0\in\bmu_{\infty}$ be a root of unity, and let $n\geq 2$ be an integer.
The symmetric group $S_n$ acts on $\widetilde{\cM}^r_{0,n+1}$ by permuting the last $n$ marked points.
We define $\cM_{n+1}(\xi_0)$ as the substack of $\bigsqcup_{h\in S_n} (\widetilde{\cM}^r_{0,n+1})^h$ made of the curves such that $\tr(\cL_0) = \xi_0^{r/\fe_0}$ (where $\fe_0$ is the order of $\mult_{x_0}(\cL)$ in $\Z/r\Z$).
\end{defn}
The polar part at $q^{1/r} = \xi_0^{-1}$ is precisely the contribution of $\cM(\xi_0)$ to Lefschetz formula.

\subsection{Expansion at $q^{1/r}=1$}

The expansion at $q^{1/r} = 1$ of the $J$-function has the form 
\begin{equation}
\label{equ:j1}
J_{(1)}(t) = 1-q+t + \sum_{\xi_0\neq 1} \mathrm{Cont}(\xi_0)(q) + \mathrm{Cont}(1)(q),
\end{equation}
where $\mathrm{Cont}(\xi_0)$ denotes the contribution of $\cM(\xi_0)$.
The function $J_{(1)}$ is an element of the space $\cK^\fake$ (see \Cref{defn:kfake}).
By definition, the only pole of $\mathrm{Cont}(\xi_0)$ is at $q^{1/r}=\xi_0^{-1}$.
Thus, in formula \ref{equ:j1}, the term $1-q+t + \sum_{\xi_1\neq 1} \mathrm{Cont}(\xi_1)(q)$ lies in $\cK^\fake_{+}$, while the term $ \mathrm{Cont}(1)(q)$ lies in $\cK^\fake_{-}$.

\begin{defn}
\label{defn:heads}
Let $C$ be an object of $\cM_{n}(1)$.
The \emph{head} of $C$ is the largest connected subcurve $C^{\head}$ such that 
\begin{itemize}
\item $x_0 \in C^{\head}$,
\item $ C^{\head}$ is $g$-stable, and $g_{|C^{\head}}=id$.
\end{itemize}
The moduli stack of heads $\cM^\head_{n_{1},n_{2}}$ is the stack parametrizing $r$-spin curves $(C,\cL,\sigma_i)$ with $n_{1}+n_{2}+1$ marked points (and section $\sigma_{i}$), together with automorphisms $\eta_i$ of $\sigma_i$, and an isomorphism $\phi:\cL\to \cL$ compatible with the spin structure.
\end{defn} 

\begin{rem}
The head of a curve $C\in \cM_{n}(1)(S)$ is always non-empty.
Indeed, since $\tr(\cL_0)=1$, the restriction of $g$ to the irreducible component containing $x_0$ is the identity.
\end{rem}
\begin{rem}
	The moduli stack of heads is the disjoint union of the $n-2$-dimensional components of $\cI \widetilde{\cM}_{0,n+1}^{r}$.
\end{rem}

Notice that on each connected component of $\cM^{\head}$, $\phi$ is the multiplication by some $r$-th root of unity $\xi$. 
Thus, the stack of heads has a natural decomposition
\begin{equation}
\cM_{n_{1},n_{2}}^\head = \bigsqcup_{\xi\in \mu_r} \cM_{n_{1},n_{2}}^\head(\xi).
\end{equation}

\begin{defn}
\label{defn:arm}
An \emph{arm}, is an object of the stack
\begin{equation}
\cM^{\arm} =  \bigsqcup_{\xi_0\neq 1}\cM(\xi_0).
\end{equation}
An arm $C$ such that $\mult_{x_0}(\cL)=0$ is called a \emph{broad arm}.
We further decompose the moduli stack of arms by taking into consideration the $g$-action on $\sigma_0^*\cL$ : 
\begin{equation}
\cM^\arm(\zeta_0)(S):=\left\{ (C,\cL,g)\in \cM^\arm(S)|\tr(\sigma_0^*\cL)=\zeta_0 \right\}.
\end{equation}
\end{defn}

Let $C$ be an object of $\cM_{n}(1)$ over a connected scheme $S$.
Then, $C$ is the union of the head and $n'$ other curves $C_1,\hdots, C_{n'}$, attached to the head at the nodes $p_1,\hdots, p_{n'}$.
Let $D_i$ be the divisor of $\widetilde{\cM}^r_{0,n}$ defined by the node $p_i$, and let $\cN_i$ be its normal bundle.
We still denote $\cN_i$ its pullback to $S$.
By definition of the head, the action of $g$ on $\cN_i$ is non trivial.
Thus, the curves $C_i$, (together with the restriction of $g$) are arms, that is, objects of $\bigsqcup_{\xi_0 \neq 1} \cM(\xi_0)$.

\begin{figure}[h]
\centering
\begin{tikzpicture}[scale=0.9]
\coordinate (a) at (-0.3,-0.3); 
\coordinate (b) at (0.8,1.3);  
\newcommand{\courbe}[1]{\draw[shift={#1}](-0.3,-0.3) ..controls +(-0.2,0.5) and +(-0.2,-0.2).. (0.8,1.3);
\draw[->, shift={#1}] (0.7,1) arc (0:340:0.2);}
\newcommand{\ligne}[1]{\draw[shift={#1}](0,-0.1) -- (0,1);}
\newcommand{\trait}[1]{\draw[shift={#1}](0,-0.1) -- (0,0.1);}

\draw[line width=0.3mm] (-0.5,0) -- (7.5,0); 
\draw (0,-0.1) -- (0,0.1);
\courbe{(2,0)}
\courbe{(1,0)}
\courbe{(4,0)}
\trait{(5,0)}
\trait{(7,0)}

\draw (0,0) node[below]{$x_0$};
\draw (1,0) node[below]{$x_1$};
\draw (2,0) node[below]{$x_2$};
\draw (3,0) node[below]{$\hdots$};
\draw (4,0) node[below]{$x_N$};
\draw (5,0) node[below]{$x_{N+1}$};
\draw (6,0) node[below]{$\hdots$};
\draw (7,0) node[below]{$x_{N+n}$};

\coordinate (origine) at (4,2);
\draw (8,0) node{Head};
\newcommand{\fleche}[1]{\draw[->] (origine) -- #1;}
\fleche{(2,1.5)}
\fleche{(3,1.5)}
\fleche{(4.5,1.5)}
\draw (origine) node[above]{Arms};

\end{tikzpicture}
\caption{Decomposition into head and arms}
\end{figure}

\begin{prop}
	\label{prop:decoheadarm}
The decomposition into head and arms yields a morphism of stacks
\begin{equation}
	\label{equ:decompositioninertie}
\bigsqcup_{\substack{ n_{1}\geq 0\\n_{2}\geq 0\\n_{1}+n_{2}\geq 2} }\cM^\head_{n_{1},n_{2}}(\xi) \times_{\cI B\bmu_r} \underbrace{\cM^{\arm}(\xi) \times_{\cI B\bmu_r} \cdots \times_{\cI B\bmu_r}\cM^{\arm}(\xi)}_{n_{1} \textrm{ times }}\to \cM(1),
\end{equation}
where the morphisms are 
\begin{subequations}
\begin{align*}
\ev_i : \cM^\head_{n_{1},n_{2}} &\to \cI B\bmu_r &\ev_0^{\vee} : \cM^{\arm} &\to \cI B\bmu_r.\\
 \end{align*}
\end{subequations}
This morphism has degree $N!(\prod_{i=1}^{n_{1}}\fe_{i})^{-1}$ onto its image, where $\fe_{i}$ is the order of the stabilizer of each node joining the head to an arm.
\end{prop}

\begin{proof} 
Let $(C^\head, C_1,\hdots, C_{n_{1}},\cL^\head,\cL^{(1)},\hdots,\cL^{(n_{1})})$ be an object of $$\cM^\head(\xi)_{n_{1},n_{2}} \times_{\cI B\bmu_r}\cM^{\arm}(\xi) \times_{\cI B\bmu_r} \hdots \times_{\cI B\bmu_r}\cM^{\arm}(\xi)$$ over a connected base scheme $S$.
Let $x_0,\hdots x_{n_{1}+n_{2}}$ be the marked points of $C^\head$, and let $y_i$ be the first marked point of $C_i$.
Each $x_i$, $y_i$ is the trivial $\bmu_{\fe_i}$-gerbe over $S$ (with $\fe_i$ the order of $\mult_{x_i}(\cL)^{\head}$), so we have a canonical isomorphism $x_i\simeq y_i$.
We define $C$ to be the gluing of the curves along this isomorphism (it exists by \cite[prop. A.1.1]{abramovichGromovWittenTheoryDeligneMumford2008}).
Moreover, we also have canonical isomorphisms $\cL^{\head}_{|x_i}\simeq \cL^{(i)}_{|y_i}$, so the line bundles glue and yield a line bundle $\cL$ over $C$.
We also have isomorphisms $(\cL^\head)^{\otimes r} \simeq (\omega_{\log})_{|C^\head}$ and $\cL^{(i)} \to (\omega_{\log})_{|C_i}$ which glue.
This is a consequence of the fact that the restriction of $\omega_{log}$ to a node or a marked point is canonically trivial, and that $\cL^{\head}_{|x_i}$ and $\cL^{(i)}_{|y_i}$ are isomorphic as maps to $B\bmu_r$.
Finally, the linearizing maps $\phi^\head,\phi_i$ glue into a global isomorphism $g^*\cL\to \cL$ because they coincide at each node. 

To compute the degree, we note that there are $n_{1}!$ choices of ordering of the nodes, and we compare the generic automorphism groups.
\end{proof}

\begin{prop}
\label{prop:facto1}
Let $\cY$ be a connected component of $\cM(1)$ such that the head of the universal curve carries $n+1$ marked points and $n_{1}$ arms. 
Let $\xi^{-1}\in\bmu_r$ be the $r$th root of unity corresponding to the morphism $\phi:\cL^\head \to \cL^{\head}$, i.e., $\tr(\sigma_0^*\cL)=\xi$.
Let $m$ be the number broad arms, $\Lambda^\head$ be the virtual class on $\cM^\head$, and $\Lambda_i$ be the virtual class on each copy of the moduli space of arms.
Then over $\cY$, the pullback by \eqref{equ:decompositioninertie} of the virtual class factorizes as follows
\begin{align*}
 \ch\left( \Tr(\Lambda)\right) &= \ch\circ\Tr\left( \Lambda^{\head}\right) \times  \prod_{i=1}^{n_{1}} \ch\left( \Tr(\Lambda_i)\right)\frac{1}{(1-s\xi)^{(N+1)n_{1}}} .
\end{align*}
Moreover, we have 
\begin{equation}
	\ch\left( \Tr\left(\Lambda^\head \right) \right) = \ch\left(\lambda_{-s\xi}R^1\pi_*\cL^\head(-E) \right)^{(N+1)}.
\end{equation}
In other words, the head contribution corresponds to the fake theory.
\end{prop}

\begin{proof}
Let $x_1,\hdots,x_{n_{1}}$ be the nodes connecting the head to the arms, and let $C^\nu = C^{\head} \sqcup \bigsqcup C_i$ be the partial normalization of the curve at these points.
Let $p:C^\nu \to C$ be the projection, $\cL^\nu = p_*p^*\cL$, and let $E^\nu$ be the divisor of broad marked points on $C^\nu$.
We write $E^\nu = E_\head \bigsqcup_{i=1}^{n_{1}} E_i$, where $E_\head$, $E_i$ are the divisors of broad marked points in $C^\head$ and $C_i$ respectively.
The divisors $E_\nu$ and $p^{*}(E)$ may differ because new broad points may arise from the normalization.
There is an exact sequence 
\begin{equation}
0\to \cL \to \cL^\nu \to \bigoplus_{i=1}^{n_{1}} \cL_{|x_i} \to 0.
\end{equation}
The pushforward $\pi_*\cL_{|x_i}$ is non-zero only if $\mult_{x_i}(\cL)=0$.
Remember that $\pi_*\cL(-E)=0$ so we have the long exact sequence 
\begin{equation}
0\to \pi_*\cL^\nu(-E) \to \bigoplus_{\fe_i=0} \pi_* \cL_{|x_i} \to R^1\pi_*\cL(-E)\to R^1\pi_*\cL^\nu(-E) \to 0
\end{equation} 
Similarly, we have a short exact sequence $0\to p_*p^*\cL(-E_\nu) \to \cL^\nu(-E) \to \bigoplus_{\fe_i =0} \cL_{|x_i}^{\oplus 2} \to 0.$ 
We get
\begin{align*}
	\MoveEqLeft[2] \lambda_{-s}(R\pi_*\cL(-E))\\
	{}&=  \lambda_{-s} (R\pi_*\cL^\nu)(-E))\bigotimes_{\fe_i=0}\lambda_{-s} (R\pi_*\cL_{|x_i})^{-1}\\
&= \lambda_{-s} (R\pi_*\cL^{\nu}(-E_\nu)) \bigotimes_{\fe_i=0}\lambda_{-s} (R\pi_*\cL_{|x_i})\\
&=  \left(\lambda_{-s} R^1\pi_*\cL^{\head}(-E_\head)\right)^{-1}\bigotimes_i\left( \lambda_{-s} R^1\pi_*\cL^{(i)}(-E_i)\right)^{-1} \bigotimes_{\fe_i=1} (1-sp_*\cL_{|x_i}).\\
\end{align*}
Since all the sequences above are exact sequences of equivariant sheaves, we may take the trace bundle and the Chern character to get the first statement 
\begin{align*}
	\ch\circ \Tr(\Lambda) &=\left(\ch\circ \Tr\left(\left(\lambda_{-s}R\pi_*\cL(-E)\right)^{-(N+1)}\right)\right)\\
			      &= \ch \circ \Tr\left( \Lambda^\head\right )\bigotimes_i \ch\circ \Tr \left( \Lambda_i \right) \bigotimes_{\fe_i=0} \frac{1}{(1-\xi s)^{(N+1)n_{1}}}.
\end{align*}
The second statement follows immediately from the assumption that $\phi$ is given by $\xi^{-1}$. 
\end{proof}

\begin{prop}
For all $t\in \left(\cK_+ \right)^{\bmu_{r}}$, we have 
\begin{equation}
J_{(1)}(t)\in \cL^{\fake}.
\end{equation}
\end{prop}

\begin{proof}
We use the decomposition into head and arms to express the terms in Lefschetz formula as an integral on the moduli space of heads.
Recall that the expansion of the $J$-function at $q^{1/r}=1$ is
\begin{equation}
J_{1}(t)=1-q +t + \sum_{\zeta\neq 1} \mathrm{Cont} (\zeta) + \mathrm{Cont} (1).
\end{equation}
Let $\widetilde{t}$ denote 
\begin{equation}
\widetilde{t} = t + \sum_{\zeta\neq 1} \mathrm{Cont} (\zeta) = \left[J_{1}(t)\right]_+ -1+q.
\end{equation}
Then we have $\widetilde{t} \in \cK^\fake_+$.
We show that $J_{1}(t) = J^\fake(\widetilde{t})$.
By the binomial formula we have
$$J^{\fake}(t+\tilde{t}) = 1-q+t+\tilde{t} + \sum_{a,n_{1},n}\frac{\phi^{a,\xi}}{\fe(a) n!n_{1}!} \scal{\frac{\phi_{a,\xi}}{1-q^{1/\fe}\cL_{0}},t,\hdots,t,\tilde{t},\hdots,\tilde{t}}_{0,n+n_{1}+1}^{\fake}.$$
Then \Cref{prop:decoheadarm} and \Cref{prop:facto1} imply that the correlator on the right and-side is equal to 
$$\int_{\cM}\ch\left( \Tr(\Lambda^{\head})\prod_{i=1}^{n_{1}}\frac{\Tr(\Lambda_{i}^{\leg})}{\Tr(1-\cL_{i+}\cL_{i-})}\right)\td(T),$$
where $\cM$ is the disjoint union of every possible product of $\cM^{\head}_{n_{1},n_{2}}$ with $n_{1}$ times $\cM^{\arm}$ as in \Cref{prop:decoheadarm}, and $\cL_{i,\pm}$ is the cotangent line at each side of a node joining an arm to the head.
Thus, by the Lefschetz formula, we have 
\begin{align*}
1-q+t+\tilde{t}+ \sum_{n,n_{1},a,\xi}\frac{\phi^{a,\xi}}{\fe(a) n!n_{1}!} \scal{\frac{\phi_{a,\xi}}{1-q^{1/\fe}\cL_{0}},t,\hdots,t,\tilde{t},\hdots,\tilde{t}}_{0,n+n_{1}+1}^{\fake} \\
=1-q+t+ \sum_{n,a,\xi} r \phi^{a,\xi} \corr{\frac{\phi_{a,\xi}}{1-q^{1/\fe}},t,\hdots,t}^{S_{n}}_{0,n+1}.
\end{align*}

\end{proof}

\subsection{Expansion at other roots of unity}

Let $\xi_0\in \bmu_{\infty}$ be a root of unity, and let $m(\xi_0)$ (or simply $m$ when $\xi_0$ is fixed) be the order of $\xi_0^r$.
Because of the $\bmu_r$-invariance, we assume that $m>1$.
The polar part of $\cJ$ at $q^{1/r} = \xi_0^{-1}$ comes from the contribution of $\cM(\xi_0)$ to the Lefschetz formula.
Since $m\neq 1$ the automorphism $g$ acts non-trivially on the connected component of $x_0$, namely, $g$ acts on this component by a rotation of order $m$.
This action allows us to decompose the curve into a spine, some legs, and a tail.
A decomposition of $\cM_{n}(\xi_0)$ as a union of products follows, and the Lefschetz formula factorizes as a product of classes called the spine/leg/tail contributions.
This allows us to recognize the expansion $J_{\xi_0}$ as a tangent vector to the Lagrangian cone of the spine CohFT.  

\begin{notation}
	We fix a root of unity $\xi_{0}$ such that $\xi_{0}^{r}$ has order $m>0$, and we write $\xi_{0} = e^{\frac{2i\pi k_{0}}{rm}}$.
By definition, $k_{0}$ is invertible modulo $m$, and we let $\nu_{0}$ denote its inverse.
\end{notation}

\begin{defn}
\label{defn:spine}
Let $C$ be an object of $\cM(\xi_0)$ over a connected scheme.
The spine of $C$ is the largest connected subcurve $C^{\spine}\subset C$, such that 
\begin{itemize}
\item $x_0\in C^{\spine}$,
\item $C^{\spine}$ is $g$-stable, and its nodes are balanced with respect to $g$,
\item $g^m_{| C^{\spine}}=\mathrm{id}$.
\end{itemize}

\end{defn}

\begin{lem}
The marked points of $C^{\spine}$ are either fixed by $g$, or have an orbit of cardinality $m$.
More precisely, the spine of a curve is isomorphic to a balanced chain of (orbifold) $\PP^1$s with standard $\Z_m$ action, with $m$-tuples of curves or marked points attached.  
Over the spine curve, the automorphism $g$ has exactly 2 smooth fixed points denoted by $x_0$ and $x_\infty$.
\end{lem}

\begin{proof}
Over the irreducible component containing the first marked point, the Riemann-Hurwitz formula implies that $g$ has exactly 2 ramification points, with maximal ramification index.
If the second ramification point is a node, then we apply the same argument to the irreducible component attached to it.  
\end{proof}

The marked points of the spine consist of one of the following
\begin{itemize}
\item $m$ tuples of permuted marked points,
\item $m$ tuples of nodes,
\item the fixed points $x_0$ and  $x_\infty$.
\end{itemize}

\begin{defn}
Let $C$ be a curve in $\cM_{n}(\xi_0)$.
If $x_\infty$ is a node of $C$, the connected component attached to $C^\spine $ at $x_\infty$ is the \emph{tail} of $C$.
The complement of the spine and the tail is a union of $\emph{legs}$.
A leg is a $m$-tuple of spin curves, cyclically permuted by $g$.

Thus the curve $C$ decomposes as the union of
\begin{itemize}
\item a spine with $mN+2$ marked points,
\item a (potentially empty) tail attached at $x_\infty$, and
\item $N$ legs, i.e. $N$ sets of $m$ cyclically permuted spin curves.
\end{itemize} 
\end{defn}
\begin{figure}[h!]
\centering
\begin{tikzpicture}[scale=1.5]
\newcommand{\trait}[1]{\draw[shift={#1}] (0,-0.1) -- (0,0.1);}
\newcommand{\traith}[1]{\draw[shift={#1}] (-0.1,0) -- (0.1,0);}
\coordinate (xi) at (-2,-0.6);
\coordinate (xhi) at (-2,0.6);
	\draw[thick,color=red] (-3,0) to[bend left=30] (1,0);
	\draw[thick,color=red] (-3,0) to[bend right=30] (1,0);

	\draw (-3,-0.1) node[below]{$x_{0}$};
	\trait{(-3,0)}
	
	\trait{(-2,0.45)}
	\trait{(-2,-0.45)}
	\draw (xi) node[below]{$x_{i}$};
	\draw (xhi) node[above]{$g(x_{i})$};

	\draw[thick] (1,0) to[out=-15, in=150] (2.5,0);

	\draw[thick] (0,0.3) to[bend left=30] (0.7,2);
	\draw[thick] (0,-0.3) to[bend right=30] (0.7,-2);
	\traith{(0.2,1.45)}
	\traith{(0.2,-1.45)}
	\draw (0.1,-1.5) node[left]{$x_{j}$};
	\draw (0.1,1.5) node[left]{$g(x_{j})$};

	\coordinate (legendetail) at (2,0.7);
	\draw (legendetail) node[above]{\textrm{Tail}};
	\draw[->] (legendetail) to (2,0.2);
	
	\coordinate (legendeleg) at (1.3,1.3);
	\draw (legendeleg) node[right]{\textrm{Leg}};
	\draw[->] (legendeleg) to (0.5,1);
	\draw[->] (legendeleg) to (0.5,-1);

	\coordinate (legendespine) at (-3,-1);
	\draw (legendespine) node[left]{\textrm{Spine}};
	\draw[->] (legendespine) to ++(0.5,0.5);
	\draw[dotted,->] (xi)+(0.1,0.2) to[bend right=40] (-1.9,0.4);
	\draw[dotted,->] (xhi)+(-0.1,-0.2) to[bend right=40] (-2.1,-0.5);
	\draw (-2,0) node{$g$};
\end{tikzpicture}
\caption{Decomposition into spine, leg and tail for $m=2$.}
\end{figure}

We know show that the decomposition of curves into spine, legs and tail leads to a decomposition of the fixed-point stack $(\widetilde{\cM}^r_{0,n})^h$ as a (union of) products of the corresponding moduli stacks.

\begin{defn}[spines]
The moduli space of \emph{spine curves} $\cM^{\spine}_{n+2}(\xi_0,\zeta_0,\zeta_\infty)$ is the stack parametrizing
\begin{itemize}
\item an $r$-spin curve $(C,\cL)$ with $mn+2$ marked points $x_{0},x_{\infty}, x_{i,k}$ ($i=1,\hdots,n$ and $k=0,\hdots,m-1$) and a section $\sigma_{i,k}$ at each marked point,
\item a balanced automorphism $g$ of $C$, with order $m$,
\item an isomorphism $\phi :g^*\cL \to \cL$ compatible with the spin structure,
\item $2$-isomorphisms $g\circ \sigma_{i,k}\simeq \sigma_{i,k+1}$,
\item $2$-isomorphisms $\eta_0 : g\circ \sigma_0 \to \sigma_{0}$ and $\eta_\infty : g\circ\sigma_\infty \to \sigma_\infty$, such that $\tr(\cL_0) = \xi_0$, $\tr(\cL_\infty)=\xi_0^{-1}$, $\tr(\sigma_0^*\cL)=\zeta_0$, $\tr(\sigma_\infty^*\cL) = \zeta_\infty$.
\end{itemize}

Notice that, on each connected component of $\cM^{\spine}$, the morphism $\phi^m : (g^m)^*\cL\to \cL$ is the multiplication by an $r$th root of unity $\xi$.

The moduli space of spines is equipped with its virtual class
\begin{equation}
\Lambda^\spine = \left( \lambda_{-s}\left(R^1\pi_*\cL(-E)\right) \right)^{\otimes N+1} \in K^0_{\Z_m}\left( \cM^\spine(\xi_0)\right)\llbracket s \rrbracket,
\end{equation}
where $E$ is the divisor of broad marked points in the universal curve.
\end{defn}

\begin{defn}[legs]
The moduli space $\cM^\leg_{n+1}$ is the stack parametrizing
\begin{itemize}
\item an $m$-tuple of $r$-spin curves $(C_i,\cL^{(i)})$ indexed by $\Z_m$,
\item isomorphisms of spin curves (with section) $g_i:(C_i,\cL^{(i)}) \to (C_{i+1},\cL^{(i+1)})$ such that the trace of $g_{m-1}\circ\hdots,g_0$ on the first cotangent line of $C_0$ is non-trivial.
\end{itemize}
We denote the composition $g_{m-1}\circ\hdots\circ g_0$ by $g^{m}$.
Let $\xi \in \bmu_r$ be an $r$th root of unity.
We denote by $\cM^{\leg}(\xi)$ the substack where $\tr_{g^m}(\sigma_0^*\cL)=\xi$.

The moduli space of legs is equipped with its virtual class
\begin{equation}
\Lambda^\leg = \left( \lambda_{-s}\left(R^1\pi_*\cL(-E)\right) \right)^{\otimes N+1}.
\end{equation}
\end{defn}

The moduli stack of legs admits a natural description in terms of arms.
\begin{lem}
\label{lem:legs}
Let $\gamma$ be the automorphism of $(\cM^\arm_{n+1})^m$ permuting the different copies of $\cM^\arm_{n+1}$.
Then, the moduli stack of legs is the fixed-point stack 
\begin{equation}
\cM^\leg_{n+1} \simeq \left( (\cM^\arm_{n+1})^m\right)^\gamma.
\end{equation}
\end{lem}

\begin{defn}[tails]
The moduli stack of tails is 
\begin{equation}
\cM^{\tail}(\xi_0,\zeta) = \bigsqcup_{\xi_0 \neq \xi_0} \cM(\xi_1,\zeta).
\end{equation}

The moduli space of tails is equipped with its virtual class
\begin{equation}
\Lambda^\tail = \left( \lambda_{-s}\left(R^1\pi_*\cL(-E)\right) \right)^{\otimes N+1}.
\end{equation}
\end{defn}

On the stack of spines $\cM^\spine(\xi_0,\zeta)$, let $\xi$ be the locally constant function corresponding to the morphism $$\phi^m:(g^m)^*\cL\to \cL.$$
\begin{lem}
\label{lem:spinexi}
We have the equality
$$\zeta^m = \xi^{-1}\xi_0^{-am},$$
where $a= \mult_{x_0}\cL$.
\end{lem}
\begin{proof}
Both sides of the equality correspond to the trace of $g^m$ on the line bundle $\sigma_0^*\cL$.
\end{proof}

We introduce a tool to glue legs to spines.
\begin{defn}
The gluing stack $G$ is defined by 
\begin{align*}
G &= \left(\left( \cI B\bmu_r \right)^m\right)^\gamma,\\
\end{align*}
where $\gamma$ is the automorphism permuting the different copies of $\cI B\bmu_r$.
\end{defn}
Objects of $G$ over a connected scheme $S$ are given by the data $(a,(L_i)_{i=0}^{m-1},(\alpha_i)_{i=0}^{m-1},$
where
\begin{itemize}
\item $a\in \Z_r$ denotes a connected component of $\cI B\bmu_r$,
\item $L_i$ is a line bundle on $S$,
\item $\alpha_i:L_i^{\otimes r}\to \cO_S$ is an isomorphism,
\item $f_i:L_i\to L_{i+1}$ is an isomorphism compatible with $\alpha_i$ (with $L_{m}=L_0$).
\end{itemize}

There are evaluation maps $\ev_{i}:\cM^{\spine}\to G$ and $\ev_{i}:\cM^{\leg} \to G$ given by the line bundles $\sigma_i^*\cL$ at a $\Z_m$-orbit of marked points.

\begin{prop}
There is a morphism of stacks 
\begin{equation}
\cM^\spine(\xi_0,\zeta_0,\zeta_\infty)_{n+2} \times_{G} \underbrace{\cM^{\leg}(\xi)\times_{G} \cdots\cM^{\leg}(\xi)}_{n \textrm{ times}} \times_{\cI B\bmu_r} \cM^\tail(\xi_0,\zeta_\infty) \to \cM(\xi_0),
\end{equation}
where $\xi = \zeta_0^m\xi_0^{a_0m}$, and the morphisms are 
\begin{itemize}
\item $\ev_i : \cM^{\spine} \to G$, $i=1,\hdots, N$,
\item $\ev_{0}^{\vee} : \cM^\leg_{n_i+1} \to G$,
\item $\ev_{0}^{\vee} : \cM^\tail_{m+1}\to \cI B\bmu_r$.
\end{itemize}
This morphism has degree $(\fe_{\infty} \prod_{i} \fe_{i})^{-1}n!m^{n}$.
\end{prop}

\begin{prop}
\label{prop:facto2}
Let $\cY$ be a connected component of $\cM^{spine}$, and $l$ (resp. $t$) be the number of broad legs (resp. tails) over $\cY$.
Then, the virtual class factorizes as 
\begin{equation}
\ch\left( \Tr \Lambda \right) = \ch\left( \Tr(\Lambda^{\spine} ) \bigotimes_i \Tr\Lambda^{\leg}_i \otimes \Tr \Lambda^{\tail} \right) \frac{1}{(1-s^m\xi)^{(N+1)l}(1-s\zeta_\infty)^{(N+1)t}}
\end{equation}
\end{prop}

\begin{proof}
Similar to the proof of \Cref{prop:facto1}.
\end{proof}

We now detail detail the contribution to the Lefschetz formula coming from the different terms in the product.

\subsection{Spine contribution}
\label{subsection:Spine contribution}
We now explicitly compute the term $\Tr(\Lambda^\spine)$ appearing in \Cref{prop:facto2}.
A spine curve $C$ is an $r$-spin curve together with a $\Z_m$-symmetry.
Quotienting by $\Z_m$, we get a curve $D$, equipped with an $rm$th root of $\omega_{\log}^{\otimes m}$.
This allows us to relate the spine contribution to the spine CohFT of \Cref{subsection:spinecohft}.

Let $C$ be an object of $\cM^{\spine}$, and let $p:C\to D=C/\Z_m$, be the quotient map, where the generator of $\Z_m$ acts by the automorphism $g$.
The algebra $p_*\cO_C$ has a decomposition $p_*\cO_C = \bigoplus_{j=0}^{m-1} T^{j}$, where $T$ is the line bundle of functions $f$ such that $g^*f=e^{2i\pi/m}f$.

\begin{lem}
Let us write $\tr(\cL_{0})=\xi_0 = e^{\frac{2i\pi k_0}{rm}}$, and let $\nu_0$ be the inverse of $k_0$ in $\Z_m$.
$T$ has non trivial multiplicities only at the two ramification points $x_0$ and $x_{\infty}$, and 
\begin{itemize}
\item $\mult_{x_0}(T)= \frac{-\nu_0}{m} $,
\item $\mult_{x_\infty}(T)= \frac{\nu_0}{m}$.
\end{itemize}
\end{lem}
\begin{proof}
	Let $\chi$ be a character of $\Z/m\Z$, and let $L_{\chi}$ be the equivariant line bundle on $C/g$ associated to $\chi$.
	Local sections of $L_{\chi}$ are functions $f$ on $C$ such that $g^{*}f = \chi(g)f$.
	This identifies $T$ with the line bundle $L_{\chi}$, with $\chi(g)=e^{\frac{2i\pi}{m}}$.
	Since $g^{-\nu_{0}}$ acts on the coarse tangent space $T_{x_{0}}C$ by $e^{\frac{2i\pi}{m}}$, we get that $\mult_{x_{0}}(L_{\chi}) = -\frac{\nu_{0}}{m}$.
\end{proof}
The couple $(\cL,\phi)$ is not always $\Z/m\Z$-equivariant sheaf.
Indeed, the morphism $\phi^{m} : (g^{m})^{*}(\cL) \to \cL$ is in general, the multiplication by an $r$th root of unity $\xi$.
Thus, for any $m$th root $\lambda$ of $\xi^{-1}$, the couple $(\cL,\lambda\phi)$ is an equivariant sheaf.
There is a natural choice for $\lambda$ given by \Cref{lem:spinexi}, ie, $\lambda = \zeta_{0} \xi_{0}^{a}$.
Since $\xi_{0}^{r}$ is a primitive $m$th root, all other roots of $\xi$ have the form $\lambda = \zeta_{0}\xi_{0}^{a+rk}$.
We define $\widetilde{\phi}= \zeta_{0} \xi_0^a \phi$, so that $(\cL,\widetilde{\phi})$ is an equivariant line bundle.
Let $\overline{\cL}$ be the line bundle over $D$ obtained by descent.
The morphism $\alpha : \cL^{\otimes r} \to \omega_{\log}$ induces an isomorphism 
\begin{equation}
\bar{\alpha} :\overline{\cL}^{\otimes rm} \simeq \omega_{log}^{m}. 
\end{equation}

\begin{prop}
	\label{prop:spine=rm}
Let us fix $\lambda = \zeta_{0}\xi_{0}^{a+rk}$. 
The data $(D,\overline{\cL},T)$ above defines a morphism 
\begin{equation*}
	f:\cM^{\spine}_{2+n}(a_0+rk,a_1,\hdots,a_\infty) \to \overline{\cM}^{rm,m}_{0,n+2}(B\Z_{m},\ua, \ub),
\end{equation*}
with $\ub = (-\nu_{0},0,\hdots, 0, \nu_{0})$, $\ua = (a_{0},ma_{1},\hdots, ma_{n}, \overline{a}_{\infty})$, where $\bar{a}_{\infty}$ is determined by the degree condition $$mn-a_0-rk-\bar{a}_\infty - \sum ma_i =0 \mod rm.$$
This morphism has degree $(\fe_{0}\fe_{\infty}\prod_{i} \fe_{i}^{m})^{-1} m^{n +1}$.

Furthermore, we have the base change
\begin{equation}
R\pi_*\cL = f^* \bigoplus_{j=0}^{m-1} R\pi_*(\overline{\cL}\otimes T^{j}),
\end{equation}
and the trace of $g$ is given by
\begin{equation}
	\Tr_g\left( R\pi_*\cL\right)= \bigoplus_{j=0}^{m-1}\zeta_{0} \xi_0^{a+rk} e^{\frac{-2i\pi j}{m}}  R\pi_*(\overline{\cL}\otimes T^{j}).
\end{equation}
In particular, the spine contribution coincides with a value of the spine CohFT (see \Cref{subsection:spinecohft}) 
$$ \ch \circ \Tr(\Lambda^\spine)  = \frac{1}{rm^{2}}f^{*}\Lambda^\spine_{\zeta_0\xi_0^{a+rk}}(\ua,\ub),$$
with $\ua = (a_0+rk,ma_1,\hdots, \bar{a}_\infty)$ and $b=(-\nu_0,0,\hdots,0,\nu_0)$.
\end{prop}
\begin{proof}	
We only prove the statement about the degree of $f$.
A general point in $\cM^{\spine}$ has an automorphism group of order $r \prod_{i}\fe_{i}\fe_{0}^{m}\fe_{\infty}$. Indeed, such an automorphism must be trivial on $C$, but may rescale $\cL$ by an $r$th root of unity, and each section comes with an automorphism group of order $\fe_{i}$.
On the other hand, an general point in $\overline{\cM}^{rm,m}_{0,n+2}$ has an automorphism group of order $rm^{2}$.
Finally, the preimage of a point in $\overline{\cM}^{rm,m}_{0,n+2}$ has $m^{n-1}$ elements because of the possible re-labellings of the marked points.
\end{proof}

\begin{rem}
	The same reasoning can also be carried at the other fixed point $x_{\infty}$, where the action of $g$ on the cotangent line is $\tr_{g}(\cL_{\infty})=\xi_{0}^{-1}$.
Let us note $\zeta_{\infty} = \tr_{g}(\sigma_{\infty^{*}\cL})$, and $\frac{a_{\infty}}{r} = \mult_{x_{\infty}}(\cL)$.
Then the roots of $\xi$ have the form $\lambda = \zeta_{\infty} \xi_{0}^{-a_{\infty}-rk}$, and the line bundle $\bar{\cL}$ obtained in this way has multiplicity $\frac{a_{\infty}}{rm}$ at $x_{\infty}$.
This justifies sending $\phi_{a,\zeta}$ to $\sum_{k=0}^{m-1} \phi_{a+rk,\zeta\xi_{0}^{-a-rk}}$ in the formula for $\Phi_{\xi_{0}}$ of \Cref{defn:morphismePhi}.
\end{rem}	
\subsection{Leg contribution}
We compute the trace bundle $\Tr(\Lambda^\leg)$, following Givental's computation in \cite{giventalPermutationequivariantQuantumKtheory2017a}.
More precisely, we want to compute the following generating function, which we call the \emph{leg contribution} :
\begin{multline}
\label{defn:legcontribution}
T_{\leg}(q)= \sum_{\substack{n\geq 2\\ \xi \in \bmu_{rm} \\ 0\leq a<r}}\frac{r\phi^{ma,\xi}}{n!}\\
\int_{\cM^{\leg}_{n+1}} \frac{\ch\left(\Tr \left(\bigotimes_{k=0}^{m-1}\ev_{k}^{*}(\phi_{a}\otimes e_{\xi^{m}}\right)  \right)\ch\left( \Tr \left(\Lambda^{\leg}\right)\right) \ch\left( \Tr \left( \bigotimes_{i=1}^{n} t(\cL_i)\right)\right)}{\ch \left(\Tr\left(\bigotimes_{k=0}^{m-1}( 1-q^{1/\fe}\cL_{0,k})\right)\right) \ch\left( \Tr \left(\lambda_{-1}\cN^\vee\right)\right)}\td(\cT),
\end{multline}

where $\cN$ is the normal bundle to the morphism $\cM^{\leg}_{n+1} \to \left(\cM_{0,n+1}^r\right)^{m}$,
$\cL_{0,k}$ are the cotangent lines at the first marked point of each component of the universal curve, and $\fe$ is the order of the stabilizer of the first marked point.
The leg contribution is a formal function of $t\in \scrI \cK_{+}$ with values in $\cK_{rm}$.

Let $(C_0,\hdots, C_{m-1})$ be an object of $\cM^{\leg}(\xi)$.
The pushforward $R^1\pi_*\cL$ is the direct sum 
\begin{equation}
\label{equation:decleg}
R^1\pi_*\cL= \bigoplus_{i\in \Z_m} E_i,
\end{equation}
where $E_i= R^1\pi_*\left(\cL(-E)\right)_{|C_i}$.
We begin with a general lemma.
\begin{lem}
\label{lem:traceadams}
Let $X$ be a smooth stack over $\C$, let $E=\bigoplus_{i\in \Z_m} E_i$ be a $\Z_m$-graded vector bundle over $X$, and let $g_i:E_i\to E_{i+1}$ be isomorphisms such that $g=g_{m-1}\circ g_{m-2}\circ\hdots \circ g_0 \in \Aut(E_0)$ has finite order.
Then we have 
\begin{equation}
\Tr_{g_\bullet} \Psi^k(E) = \left\{ \begin{matrix}
m\Psi^m\left (\Tr_g \Psi^{k/m}E_0 \right)& \textrm{ if } m|k\\
0 &\textrm{ otherwise.}
\end{matrix}\right.
\end{equation}
\end{lem}
\begin{proof}
If $g=\mathrm{id}$, then this is the situation in \cite{giventalPermutationequivariantQuantumKtheory2017a}.
Let us recall the argument.
There is an isomorphism of $\Z_m$-bundles $E \simeq E_0\otimes \cO_{X}[\Z_r]$, and we compute that $\Psi^k(\cO_{X}[\Z_r])=\cO_{X}^m$ is $m|k$, and $\Psi^k(\cO_{X}[\Z_r])=0$ otherwise.

If $g\neq \mathrm{id}$, we can decompose each $E_i$ into the sum of eigenspaces for $g$.
These eigenspaces are preserved by the $g_i$, so we may assume that $g=\lambda \mathrm{id}$ for some $\lambda \in \C^*$.
Then we have 
\begin{equation}
\Tr_{g_\bullet}\left( \Psi^k(E)\right) = \left\{\begin{matrix}
m\lambda^{k/m}\Psi^k(E_0) &\textrm{ if } m|k,\\
0 & \textrm{ otherwise.}
\end{matrix}\right.
\end{equation}
\end{proof}

\begin{cor}
We extend the Adams operations to $K^0(X)\llbracket s \rrbracket$ by setting $$\Psi^m(s)=s^m.$$
Then, with the same notations as in the previous lemma, we have  
\begin{equation}
\Tr\left( \lambda_{-s} E \right)= \Psi^m \left( \Tr\left( \lambda_{-s}E_0\right)\right).
\end{equation}
\end{cor}
\begin{proof}
Recall that $\lambda_{-s}(E)= \exp\left( -\sum_{k\geq 1} \frac{s^k\Psi^m(E)}{k}\right)$.
Thus, using \Cref{lem:traceadams} we get   
\begin{align*}
\Tr\left( \lambda_{-s}E\right)&= \exp\left(- \sum_k \frac{s^k\Tr\Psi^k\left(\bigoplus E_i\right)}{k}\right)\\
&=  \exp\left(- \sum_k \frac{s^{km}\Psi^m \Tr\Psi^k\left( E_0\right)}{k}\right)\\
&= \Psi^m \left( \Tr_{g}\lambda_{-s}E_0\right).
\end{align*}
\end{proof}

\begin{prop}[\cite{giventalPermutationequivariantQuantumKtheory2017a} lemma p.5]
Let $X$ be a proper smooth Deligne-Mumford stack over $\C$, and let $\pi:X\to \Spec( \C)$ be the projection.
Let $E,T $ be vector bundles over $X$, and let $E_i,T_i$ be the pullbacks $E_i=p_i^*E$, and $T_i=p_i^*T$,  where $p_i:X^m\to X$ is the i-th projection.
Let $\phi:X^m \to X^m$ be the cyclic permutation of factors.
We choose  finite order isomorphisms $g_i:E\to E$ and $f_i:T\to T$, which induce isomorphisms $g_i:\phi^*E_i\to E_{i+1}$ and  $f_i:\phi^*T_i \to T_{i+1}$.
We equip $E_\bullet=\bigoplus_{i=0}^{m-1} E_i$ and $T_\bullet=\bigotimes_{i=1}^{m-1} T_i$ with the induced $\Z_m$ equivariant structure. 
Finally, let $g:E\to E$ (resp. $f:T\to T$) be the composition $g=g_{m-1}\circ \cdots \circ g_0$ (resp. $f=f_{m-1}\circ \cdots \circ f_0$).

Then we have
\begin{equation}
\tr H^*\left(X^m; \left( \Psi^k\left(\frac{ E_\bullet}{k}\right)\right)\otimes T_\bullet \right) 
= \left\{ \begin{matrix}
\Psi^m\left(H^*\left(X;\Tr_g\frac{\Psi^{k/m}E}{k/m}\otimes \Tr_f T\right)\right) &\textrm{ if } m|k,\\
0&\textrm{ otherwise.}
\end{matrix}\right.
\end{equation}

\end{prop}
\begin{proof}
Use Lefschetz formula to compute the left hand-side, and Adams--Riemann--Roch for the right hand-side, combined with \Cref{lem:traceadams}.
\end{proof}

We apply the previous results to the space $\cM^\leg \simeq \left( (\cM^{\arm})^m\right)^\phi$ (see \Cref{lem:legs}).
Let $\pr$ denote the projection to the first factor.
\begin{cor}
We have 
\begin{equation}
\tr\left( H^*(\Lambda^\leg)\right)= \Psi^m\left( H^*(\Tr(\Lambda^\arm)\right) \in \C\llbracket s  \rrbracket.
\end{equation}
More generally, the leg contribution is given by
$$\Psi^m\Phi_{0}\left( \left[J_{(1)}(t)\right]_+\right),$$
where $[\cdot]_+$ denotes the projection to $\cK^\fake_+$, parallel to $\cK^\fake_-$.
\end{cor}
\begin{proof}
	This follows from the previous computations and the equality 
	\begin{equation}
	\Psi^{m}\left(\phi_{a,\xi}\right) = \sum_{\zeta^{m} = \xi} \phi_{a,\xi}.
	\end{equation}
\end{proof}
\subsection{$J_{(\xi_{0})}$ is tangent to $L^{\spine}$}
We now explain how the previous results imply that $J_{\xi_0}$ is a tangent vector to the cone $L^\spine$ of the spine CohFT.
Recall that its state space is $K^0(IB\bmu_{rm})\otimes H^*(\cI B\Z_m,\C) \otimes \C\llbracket s \rrbracket$, with the orbifold pairing.
The factorization of the virtual class means that we can view the leg and tail contributions as inputs of the spine CohFT.

Let us decompose $H^*(\cI B\Z_m,\C) = \bigoplus_{d\in\Z_m} \C \cdot [d]$.
We refer to the subspaces $K^0(\cI B\bmu_{rm})_\C \otimes [d] $ as sectors of the state space.
Recall (\Cref{defn:morphismePhi}) that the two embeddings $\Phi_0,\Phi_{\xi_0}$ of $\cK$ into $\cK_{rm}$ correspond to the sectors $0$ and $\nu_0 = \frac{1}{k_0}\mod m$ respectively (with $\xi_0=\exp({\frac{2i\pi k_0}{rm}})$).

\begin{prop}
\label{prop:image de J tangent}
For all $t\in \left( \cK_{+}\right) ^{\bmu_{r}}$, the image of $J(t)(q^{1/rm})$ by $\Phi_{\xi_0}$ is a tangent vector to the spine cone :
\begin{equation}
\Phi_{\xi_0}\left( J(t)\right) (q^{1/rm})\in \cT L^\spine.
\end{equation}
The tangency point is $J^\spine(T)$, with $$T=\Psi^m\Phi_0\left(\left[J_{(1)}(t)\right]_{+} \right).$$
\end{prop}

Before proving this proposition, we need a preliminary lemma.
\begin{lem}
Let $C$ be an object of $\cM_{mn+2}^\spine(\xi_0,\zeta_0,\zeta_\infty)$, and let $\ua \in \Z_r^n$ be the multi-index of multiplicities of $\cL$ at each $\Z_m$-orbit of marked points.
Then we have
\begin{equation}
\zeta_\infty = \zeta_{0}\xi_0^{m\sum a_i -mn}
\end{equation}
\end{lem}

\begin{proof}
Let us assume that $\xi_0=e^{\frac{2 i\pi k_0}{rm}}$ is a primitive $rm$th root of unity, and let $\nu$ be the inverse of $k_0$ modulo $rm$.
We may also assume that $(\cL,\phi)$ is a $\Z_m$-equivariant bundle, and that $C$ is smooth, that is, $C$ is isomorphic to a stacky $\PP^1$.
Let $\bar{L}$ be the bundle over $D=C/\Z_m$ obtained by descent, and let $\bar{a}_0,\bar{a}_\infty,\bar{a}_i$ be its multiplicities at the marked points.

Let $p:U\to C$ be the $\Z_r$-cover ramified over $x_0$ and $x_\infty$.
Then, there exists a lift $\bar{g}$ of $g$ such that
\begin{itemize}
\item $\bar{g}$ has order $rm$, and its trace on the tangent space at $x_0$ and $x_\infty$ is $\xi_0$ and $\xi_0^{-1}$ respectively,
\item $p^*\cL$ is $\Z_{rm}$-equivariant, and descends to $\overline{\cL}$ on the quotient $U/\Z_{rm}\simeq D$.
\end{itemize}
Then $\zeta_0$ and $\zeta_\infty$ are the trace of $\bar{g}$ on $p^*\cL$ at $x_0$ and $x_\infty$ respectively.
Thus, we have that 
\begin{align*}
\zeta_0^{-\nu} &= e^{\frac{2i\pi \bar{a}_0}{rm}},\\
\zeta_\infty^{\nu} &= e^{\frac{2i\pi \bar{a}_\infty}{rm}}.
\end{align*}
Finally, we use the fact that $\bar{a}_0 + \bar{a}_\infty = \sum_i m a_i -mn$ to obtain
\begin{align*}
\zeta_{\infty} &= \xi_0^{\bar{a}_\infty}\\
&=\zeta_{0}\xi_0^{m\sum a_i -mn}.
\end{align*}
The proof is similar for a general $\xi_0$. 
\end{proof}

\begin{proof}[Proof of \Cref{prop:image de J tangent}]
Let $\delta t$ be the tail contribution, namely
\begin{equation}
	\delta t = \sum_{\substack{a\in \Z_r\\ \xi \in \bmu_r}}\sum_{n\geq 2} \frac{r\phi^{a,\xi}}{ n!} \int_{\cM^\tail(\xi_0)} \frac{\ch\left( \Tr\left(\Lambda^\tail\otimes \ev_0^*(\phi_{a}\otimes e_\xi) \bigotimes t(\cL_i)\right)\right)}{\left( 1-q^{1/m\fe(a)}\Tr(\cL_0)\right) \ch\left( \Tr(\lambda_{-1}\cN^\vee)\right)} \td(\cT),
\end{equation}
where $\cN$ is the normal bundle to the morphism $\cM_{n+1}^\tail \to \overline{\cM}^{r}_{0,n+1}$, and $\cT$ is the tangent bundle.
Then we have 
$$\delta t(q^m) = \left[J(t)(q^{1/r}\xi_0^{-1})\right]_+.$$ 
We claim that
\begin{equation}
\Phi_{\xi_0}\left( J(t)(q^{1/rm})\right) = \Phi_{\xi_0}(\delta t) \\+\sum_{\substack{a\in \Z_{rm} \xi \in \bmu_{rm}}} \sum_{n\geq 1} \frac{\phi^{a,\xi}\otimes [\nu_0]}{m\fe(a) n!} \scal{\frac{\phi_{a,\xi} \otimes [-\nu_0] }{1-q^{1/\fe(a)m}},T,\hdots,T,\Phi_{\xi_0}(\delta t)}^{\spine}_{0,n+2},
\end{equation}
with $T = \Psi^m \left(\Phi_0\left( \left[J_{(1)}(t)\right]_+\right)\right)$, and $\fe(a)$ is the order of $a$ in $\Z_r$.
This claim follows from the factorization of the virtual class, \Cref{lem:invariancespine}, \Cref{prop:spine=rm}, and the fact that the conjugacy class of the partition $N=\sum_{l\geq  1} l n_{l}$ in $S_{N}$ has $\frac{N!}{\prod_{l\geq  1} l^{n_{l}}n_{l}!}$ elements.
Thus, $\Phi_{\nu_0}\left( J(t)\right)$ is a tangent vector to $L^\spine$ at $\Psi^m(\Phi_0 J(t))$.
\end{proof}

Let $\cT^{\nu_0}$ denote the tangent space to $L^{\spine}$ at $\Psi^m\Phi_0(J_{1}(t))$ intersected with the $\nu_0$ sector. 
Then we have
\begin{equation}
\cT^{\nu_0} = \Box_{\xi_{0}} \cT L^\un,
\end{equation}
where the second tangent space is computed at $\Box_0^{-1}\Psi^m\Phi_0(J_{1}(t)) = \Delta^{-1}\Phi_0(J_{1}(t))$ (see \Cref{prop:twistedfakepsim}), and $L^\un$ is the untwisted cone in $\cK_{rm}$.

\begin{prop}
We have 
$$ \Phi_{\xi_0}J(t)(q^{1/rm}) \in \Box_{\xi_{0}} \Delta^{-1} \cT_{\Phi_0J_{(1)}(t)}L^\fake$$
\end{prop}

\begin{proof}
This is a consequence of \Cref{prop:psim(delta)=box0}.
\end{proof}

\subsection{Reconstruction}
\label{subsection:reconstruction}
So far, we proved that the values of the $J$-function satisfy the 3 conditions in \Cref{thm:adelic_char}.
We now explain how these properties allow to reconstruct the $J$-function.
We follow the proof of \cite[prop. 4]{giventalHirzebruchRiemannRoch2011}.
Recall that the ground $\lambda$-ring $R$ is supposed to carry a Hausdorff $\scrI$-adic topology such that $\Psi^m\left( \scrI\right) \subset \scrI^{m}$.

Let $f$ be a $\bmu_r$-invariant element of $\cK$ satisfying the conditions of \Cref{thm:adelic_char}.
We write $f=1-q +t + f_-$, where $1-q+t\in \cK_+$, and $f_-\in \cK_-$. 
By assumption, we have that $f=1-q \mod \scrI$.
Notice that the last 2 conditions of \Cref{thm:adelic_char} are stable by base change.
In particular, if $f$ is an element satisfying those conditions, then so does the image of $f$ modulo $\scrI^n$.
We will show by recursion that for all $n\in \N$, $f$ is a value of the $J$-function modulo $\scrI^n$.

For $n=1$, we assumed that $f=1-q \mod \scrI$, which is a value of the $J$-function.

Now, suppose that $f=J(t) \mod \scrI^n$ for some $n\geq 1$, and let us show that $f=J(t) \mod \scrI^{n+1}$.
We just need to check that $\tilde{t} = \left[ f_{(1)}\right]_+-1+q-t$ is the arm contribution modulo $\scrI^{n+1}$, \emph{i.e.}, that we have $\tilde{t} = [J_{(1)}(t)]_+ -1+q-t$.
By definition the arm contribution is the sum of the polar parts of $J(t)$ at all the non-trivial roots of unity, so we need to show that for all $\xi_0\in \bmu_\infty$, the polar part of $f$ at $q^{1/r} = \xi_0^{-1}$ matches the polar part of $J(t)$. 
Let us begin with the $r$th roots of unity.
First, notice that $[f_{(1)}]_- \mod \scrI^{n+1}$ is determined by $[f_{(1)}]_+ \mod \scrI^{n}$ by the formula $f_{(1)} = J^\fake( [f_{(1)}]_+-1+q)$.
The $\bmu_r$-invariance of $f$ implies that the polar part at $q^{1/r} = \xi^{-1}\in \bmu_r$ is exactly $\xi\cdot \left[f_{(1)}\right]_-$.
Since $J(t)$ is also $\bmu_r$-invariant, its polar part at the $r$th roots of unity coincide with that of $f$.
Now, suppose that $\xi_0$ is some root of unity such that the order of $\xi_0^r$ is $m \geq 2$.
Then, the polar part of $f$ at $q^{1/r}=\xi_0^{-1}$ is determined by $\left[\Phi_{\xi_0} f \right]_- \mod \scrI^{n+1}$, which only depends on the tail contribution modulo $\scrI^{n}$, and on $[f_{(1)}]_+$ modulo $\scrI^{\lfloor \frac{n+1}{m}\rfloor +1}$ because of the third condition of the adelic characterization.
Thus, the induction hypothesis allows us to conclude that the polar parts of $f$ and $J(t)$ coincide modulo $\scrI^{n+1}$ at any root of unity.

Finally, we conclude that $f=J(t)$, which concludes the proof of \Cref{thm:adelic_char}.

\section{$I$-function and difference equation}
In this section, we use the adelic characterization to give a simpler description of the image of the $J$-function using ``untwisted" invariants.
Then, we use this description to find a specific value of the $J$-function, following the method given by Coates, Corti, Iritani and Tseng in \cite{coatesComputingGenusZeroTwisted2009}
In this section we choose $R$ to be the $\lambda$-ring $\C[X]$, with Adam's operations $\Psi^k(X)=X^k$.

\subsection{FJRW invariants from untwisted invariants}
\label{subsection:untwisted}
\begin{defn}
Let $n=k_1+\cdots + k_s$ be a partition of $n$, and let $H\subset S_n$ be the subgroup $H=S_{k_1}\times \cdots S_{k_n}$.
For a sequence $t^{(1)},\hdots, t^{(s)}$ of elements of $\cK_+$, the cohomology group
\begin{equation}
\left[ t^{(1)}(\cL_{1,1}),\hdots, t^{(1)}(\cL_{1,k_1});\hdots,t^{(s)}(\cL_{s,k_s})\right]^\un_{n}=H^*\left( \widetilde{\cM}^r_{0,n} ; \bigotimes_{k=1}^s \bigotimes_{l=1}^{k_s}t^{(k)}(\cL_{k,i})\right)
\end{equation}
is an $H$-module.

For any elements $\nu_1,\hdots, \nu_s \in R$, we define the untwisted invariants by
\begin{multline}
\scal{t^{(1)}(\cL_{1,1})\otimes \nu_1, \hdots ,t^{(s)}(\cL_{s,k_s})\otimes \nu_s}_{0,n}^{\un,H}  =\\ \frac{1}{\prod_i k_i!}\sum_{h\in H}\tr_h \left[ t^{(1)}(\cL_{1,1}),\hdots,t^{(s)}(\cL_{s,k_s})\right]^\un_{n} \prod_{i=1}^s \prod_{r=1}^\infty \Psi^r(\nu_i)^{l_r(h)}.
\end{multline}
The associated $J$-function is $$J^{K,\un}(t) = 1-q+t + \sum_{a,\xi} \phi^{a,\xi} \scal{\frac{\phi_{a,\xi}}{1-q^{1/\fe(a)}\cL_0},t(\cL_1),\hdots,t(\cL_n)}_{n+1}^{S,\un},$$
and the image of the $J$-function is the subvariety $$L^{K}_{\un}\subset \cK.$$
\end{defn}

\begin{thm}[Adelic characterization]
Let $f=1-q+t+f_-$ be a $\bmu_r$-invariant element of $\cK$ such that $t + f_{-}\in \scrI \cK$.
Then, $f$ belongs to $L^{K}_{\un}$ if and only if 
\begin{itemize}
\item the poles of $f$ are at $q=0$, $q=\infty$, and the roots of unity,
\item $f_{(1)}\in L^{H}_{\un}$,
\item $\Phi_{\xi_0} (f)(q^{\frac{1}{rm}}) \in T$, where $T$ is the tangent space to the untwisted cohomological cone at the point $J^{H,\un}(t)_{(1)}$.
\end{itemize}
\end{thm}
\begin{proof}
The proof of the previous adelic characterization carries verbatim, up to a change of inner product which occurs because we replace $\Lambda$ by $1$.
\end{proof}

\begin{thm}
\label{thm:LdeltaLun}
The invariant points of $L$ coincide with the invariant points of $\Delta L^{K}_{\un}$ :
\begin{equation}
(L_{\FJRW}^K)^{\bmu_r} = \left( \Delta L^{K}_{\un}\right)^{\bmu_r},
\end{equation}
where $\Delta$ is the operator of the fake theory.
\end{thm}

\begin{lem}
\label{lem:boxdelta}
Let $f$ be a point of $\cK$, and $\xi_0$ a root of unity such that $\xi_0^r$ has order $m$.
Then we have
\begin{equation}
	\Phi_{\xi_0}\left(  \Delta f \right) (q^{1/rm})= \Box_{\xi_0} \left( \Phi_{\xi_0}f(q^{1/m})\right).
\end{equation} 
\end{lem}
\begin{proof}
This is a direct computation.
\end{proof}

\begin{proof}[Proof of \Cref{thm:LdeltaLun}]
We check the criteria of the adelic characterization.
First, a direct computation shows that $\Delta$ sends $\bmu_r$-invariant points to $\bmu_r$-invariant points.
The first 2 items of the adelic characterization are obviously satisfied, and the last one follows from \Cref{lem:boxdelta}.
\end{proof}

\subsection{The $I$-function}

We begin by calculating a point of the untwisted permutation-equivariant cone.
\begin{prop}
	We have 
	\begin{equation}
		I^\un :=J^{S,\un}\left(x\phi_{2,0}\right) = (1-q) \phi_{1,0} + (1-q)\sum_{n\geq 1 } \frac{x^n}{\prod_{k=1}^n(1-q^k)}\phi_{n+1,0}.
\end{equation}
\end{prop}
\begin{proof}
Recall \cite{giventalPermutationEquivariantQuantumKtheory} that the $J$-function of a point is 
\begin{equation}
J_{\pt}(x)= 1-q + (1-q)\sum_{n\geq 1 } \frac{x^n}{\prod_{k=1}^n(1-q^k)}.
\end{equation}
In our case, the untwisted invariants coincide with the invariants of a point

\begin{equation}
	\scal{\frac{\phi^{n+1,0}}{1-q\cL_{0}^{\fe(n+1)}},x\phi_{2,0},\hdots, x\phi_{2,0}}^{\un,S_{n}}_{0,n+1} = \scal{\frac{1}{1-q\cL}, x,\hdots, x}^{\mathrm{pt,S_{n}}}_{0,n+1}.
\end{equation}
So we only need to prove that 
\begin{equation}
	\scal{\frac{\phi^{n+1,0}}{1-q^{\frac{1}{\fe(n+1)}}\cL_{0}}g,x\phi_{2,0},\hdots, x\phi_{2,0}}^{\un,S_{n}}_{0,n+1}=\scal{\frac{\phi^{n+1}}{1-q\cL_{0}^{\fe(n+1)}},x\phi_{2,0},\hdots, x\phi_{2,0}}^{\un,S_{n}}_{0,n+1} . 
\end{equation}
Both terms are rational functions of $q^{\frac{1}{\fe(n+1)}}$ so we may check that they are equal by expansion as formal power series of $q^{\frac{1}{\fe(n+1)}}$.
Then the result follows from \Cref{rem:vanishing}.
\end{proof}

We now follow the computation in \cite{coatesComputingGenusZeroTwisted2009},\cite{chiodoLandauGinzburgCalabiYauCorrespondence2010}.
We define $w_\xi(z)= \sum_{d\geq 0} w_{\xi,d} \frac{z^d}{d!}$ with $$w_{\xi,d} = (N+1)\sum_{k\geq 1} \frac{\xi^{k}s^kk^d}{k}.$$
Notice that $\exp(-w_\xi(z))= \left( 1-\xi sq\right)^{N+1}.$
We also define the functions $G_{y,\xi}$ by
\begin{align*}
G_{y,\xi}(x,z)= \sum_{k,l} w_{\xi,k+l-1}\frac{B_l(y) x^kz^{l-1}}{k!l!}.
\end{align*}
These functions satisfy two equations :
\begin{align*}
G_{y,\xi}(x,z) &= G_{0,\xi}(x+yz,z),\\
G_{0,\xi}(x+z,z)&= G_{0,\xi}(x,z) + w_\xi(x).
\end{align*}
Let $\nabla$ be the vector field $\nabla = \frac{x}{r}\partial_x$.
\begin{prop}
The function $\exp\left(-G_{1/r}(z\nabla,z)\right)I^\un$ is a value of the untwisted permutation-equivariant $J$-function $J^{S,\un}$.
\end{prop}
\begin{proof}
We check the conditions of the adelic characterization.
We compute that 
\begin{equation}
\exp\left(-G_{1/r}(z\nabla,z)\right) x^n \phi_{n+1,\xi}= \exp\left(-(N+1)\sum_{k\geq 1} \frac{\xi^{k}s^{k}}{k} \frac{q^{k\frac{n+1}{r}}}{q^{k}-1}\right) x^n \phi_{n+1,\xi}.
\end{equation}
This shows that $\exp\left(-G_{1/r}(z\nabla,z)\right)I^\un$ is $\bmu_r$-invariant, and has poles at $0,\infty$, and the roots of unity.

By a theorem  Coates--Corti--Iritani--Tseng \cite[thm. 4.6]{coatesComputingGenusZeroTwisted2009}, the operator $\exp\left(-G_{1/r}(z\nabla,z)\right)$ preserves the cone $L^{H,\un}$.
Thus, the second condition is also satisfied.
The same argument also applies to show that the third conditions holds.
\end{proof}

\begin{thm}
\label{thm:fonctionI}
	The function
$$I^K_{\FJRW}=(1-q)\sum_{\xi\in\bmu_r} \sum_{n\geq 0} \frac{\prod_{0\leq k <\lfloor n/r\rfloor} \left(1-\xi sq^{ \left\{\frac{n}{r} \right\}+\frac{1}{r}+k}\right)^{N+1}x^n}{\prod_{k=1}^n (1-q^k)} \phi_{n+1,\xi}$$
is a point of the image $L^K_{\FJRW}$ of the $J$-function.
\end{thm}
\begin{proof}
We compute $I^K_{\FJRW}:=\Delta\exp\left(-G_{1/r}(z\nabla,z)\right)I^\un$.
For $n+1\notin r\Z$ we have $\Delta \phi_{n+1,\xi} = \exp\left( G_{0,\xi}(\frac{n+1}{r}z,z)\right)\phi_{n+1,\xi}$, so
\begin{align*}
	\MoveEqLeft[3] \Delta \exp\left(-G_{1/r,\xi}(z\nabla,z)\right) x^{n} \phi_{n+1,\xi}\\ 
	{}&= \exp\left( G_{0,\xi}\left(\left\{\frac{n+1}{r}\right\}z,z\right) -G_{0,\xi}\left(\frac{n+1}{r}z,z\right)\right)x^{n} \phi_{n+1,\xi}\\
	  &= \prod_{k=0}^{\lfloor n/r \rfloor-1} \exp\left(-w_\xi\left( \left\{\frac{n+1}{r}\right\}z+kz\right)\right)x^{n} \phi_{n+1,\xi} \\
	  &=  \prod_{k=0}^{\lfloor n/r \rfloor-1} \left( 1-s \xi q^{\left\{\frac{n}{r}\right\} + \frac{1}{r}+k}\right)^{N+1}x^{n} \phi_{n+1,\xi}.
\end{align*}

For $n+1\in r\Z$, we have 
\begin{align*}
	\Delta \phi_{0,\xi} &= \exp\left( G_{0,\xi}(z,z)\right).\\
\end{align*}
Thus we have
\begin{align*}
\Delta \exp\left(-G_{1/r,\xi}(z\nabla,z)\right) x^{n} \phi_{0}\otimes e_\xi &= \exp\left( G_{0,\xi}\left(z,z\right) -G_{0,\xi}\left(\frac{n+1}{r}z,z\right)\right)x^{n} \phi_{n+1}\otimes e_\xi \\
&=\prod_{k=0}^{\lfloor n/r \rfloor-1} \exp\left(-w_\xi\left( z+kz\right)\right)x^{n} \phi_{n+1}\otimes e_\xi \\
&= \prod_{k=0}^{\lfloor n/r \rfloor-1} \left( 1-s\xi q^{1+k}\right)^{N+1}.
\end{align*}
\end{proof}

\subsection{Difference equation}
In this section we take the limit $s=1$ and we give the difference equation satisfied by the $I$-function.
Up to a change of variable, this equation coincides with that satisfied by hypersurfaces of degree $r$ in $\PP^{N}$.

We expand the $I$-function with respect to the basis $\phi_{a,\xi}$ $$I_{\FJRW}(x,q) = \sum_{a=0}^{r-1}\sum_{\xi \in \mu_r} x^{a}I_{a,\xi}(x,q) \phi_{a+1,\xi},$$
and we introduce the modification $$\tilde{I}_{a,\xi}(x,q^{-1}) := e_{q,q^{\frac{a+1}{r}}\xi^{-1}}(x)I_{a,\xi}(x^{1/r},q^{-1}),$$
where $e_{q,\lambda}(x)$ is the $q$-character satisfying the equation $q^{x\partial_{x}}e_{q,\lambda} = \lambda e_{q,\lambda}$.

\begin{thm}
	\label{thm:qdiffeq}
	The functions $\tilde{I}_{a,\xi}(x,q^{-1})$ are solutions to the $q$-difference equation
	\begin{equation}
	\label{equ:qdiffFJRW}	
	\left[ \prod_{k=1}^{r}\left( 1-q^{rx\partial_{x}-k}\right) +x(-1)^{r+N}q^{\frac{r(r-1)}{2} + (r^{2}-N-1)x\partial_{x}} (1-q^{x\partial_{x}})^{N+1}\right]I = 0.
	\end{equation}
\end{thm}
\begin{proof}
	We have
	$$I_{a,\xi}(x^{1/r},q^{-1})= (1-q^{-1}) \sum_{d \geq 0} \frac{\prod_{0\leq k<d}\left( 1-\xi q^{-\frac{a+1}{r}-k}\right)^{N+1}}{\prod_{k=1}^{rd+a}(1-q^{-k})}x^{d}.$$
This functions satisfies the equation
\begin{equation}
	\left[\prod_{k=0}^{r-1}\left( 1-q^{-r x\partial_{x}-a+k}\right)-x \left( 1-\xi q^{-\frac{a+1}{r}-x\partial_{x}}\right)^{N+1} \right]I_{a,\xi}(x,q^{-1})=0,
\end{equation}
which is equivalent to 
\begin{multline}
	\left[\prod_{k=0}^{r-1}\left( 1-q^{r x\partial_{x}+a-k}\right) +(-1)^{r+N} x q^{r^{2}+r^{2}x\partial_{x}+ra -\frac{r(r-1)}{2}} \left(\xi q^{-\frac{a+1}{r}-x\partial_{x}}\right)^{N+1} \right. \\ \left. \left( 1-\xi^{-1} q^{\frac{a+1}{r}+x\partial_{x}}\right)^{N+1} \right]I_{a,\xi}(x,q^{-1})=0.
\end{multline}
Multiplying be the exponential $e_{\xi^{-1} q^{(a+1)/r}}$, we obtain the desired equation
\begin{equation}
	\left[\prod_{k=1}^{r}\left( 1-q^{r x\partial_{x}-k}\right) +(-1)^{r+N} x q^{(r^{2}-N-1)x\partial_{x}+\frac{r(r-1)}{2}} \left( 1-q^{x\partial_{x}}\right)^{N+1} \right]\tilde{I}_{a,\xi}(x,q^{-1})=0.
\end{equation}
\end{proof}
\begin{rem}
	The $r^{2}$ components $\tilde{I}_{a,\xi}$ are exactly the solutions computed by Wen in \cite{wenDifferenceEquationQuintic2022}.
	The $I$-function also coincides with that of Aleshkin--Liu \cite{aleshkinWallcrossingKtheoreticQuasimap2022}.
\end{rem}

\begin{prop}[\cite{wenDifferenceEquationQuintic2022}]
	If $N+1 \leq r^{2}$ then the functions $\tilde{I}_{a,\xi}$ form a fundamental system of solutions to the difference equation \eqref{equ:qdiffFJRW}.
\end{prop}

Finally, let us describe the relation of \Cref{thm:fonctionI} to the quantum K-theory of hypersurfaces.
The $I$-function of a hypersurface $X\subset \PP^{N}$ of degree $r$ is known to be \cite{giventalPermutationEquivariantQuantumKtheoryb,tonitaTwistedKtheoreticGromov2018}
\begin{equation}
	I_{X}^{K} = (1-q) \sum_{d\geq 0} Q^{d} \frac{\prod_{k=1}^{dr}(1-P^{l}q^{r})}{\prod_{k=1}^{d}(1-Pq^{r})^{N+1}},
\end{equation}
which is a solution to the equation
\begin{equation}
\label{equ:qdiffhypersurface}
\left[ (1-Q\partial_{Q})^{N+1} - Q \prod_{k=1}^{r}(1-q^{k+rQ\partial_{Q}})\right]I=0.
\end{equation}

After the change of variable $Q = x^{-1}$, \eqref{equ:qdiffhypersurface} coincides with \eqref{equ:qdiffFJRW}. 

\bibliographystyle{alpha}
\bibliography{Biblio}
\bigskip

{\sc Sorbonne Université and Université Paris Cité, CNRS, IMJ-PRG, F-75005 Paris, France}

{\sc Academia Sinica, Taipei, Taiwan}

E-mail address: cazaux@as.edu.tw

\end{document}